\documentclass{amsart}

\usepackage{amsmath,amsthm,amssymb,amscd,mathtools}
\usepackage{graphicx} 
\usepackage{amsfonts}
\usepackage{rotating}
\usepackage{euscript}
\usepackage{pst-node}
\usepackage{epsfig,verbatim}
\usepackage{caption}
\usepackage{enumerate}
\usepackage{amsthm}
\def\be{\begin{equation}}
\def\ee{\end{equation}}

\def\C{{\mathbb C}} 
\def\f{\EuScript}
\def\N{{\mathbb N}} 
\def\P{{\mathbb P}}

\def\Q{{\mathbb Q}}

\def\e{\eqref}
\def\phi{{\varphi}}
\def\v{{\varepsilon}} 
\def\tt{\widetilde}
\def\deg{{\rm deg\,}}

\def\Ker{{\rm Ker\,}} 
\def\Gal{{\rm Gal\,}}

\def\GCD{{\rm GCD }}
\def\LCM{{\rm LCM }}

\def\bp{\begin{proposition}}
\def\ep{\end{proposition}}

\def\bt{\begin{theorem}}
\def\et{\end{theorem}}
\def\br{\begin{remark}}
\def\er{\end{remark}}
\def\be{\begin{equation}}
\def\l{\label}
 
\def\e{\eqref}

\def\ee{\end{equation}}
\def\eee{\end{equation*}}
\def\bl{\begin{lemma}}
\def\el{\end{lemma}}
\def\bc{\begin{corollary}}
\def\ec{\end{corollary}}
\def\pr{\noindent{\it Proof. }}

\def\bd{\begin{definition}}
\def\ed{\end{definition}}
\def\bpr{\begin{problem}}
\def\epr{\end{problem}}
\def\bco{\begin{conjecture}}
\def\eco{\end{conjecture}}
\def\t{\widetilde}

\newtheorem{theorem}{Theorem}[section]
\newtheorem{lemma}[theorem]{Lemma}
\newtheorem{definition}[theorem]{Definition}
\newtheorem{corollary}[theorem]{Corollary}
\newtheorem{proposition}[theorem]{Proposition}
\newtheorem{remark}[theorem]{Remark}
\newtheorem{problem}[theorem]{Problem}
\newtheorem{conjecture}[theorem]{Conjecture}

\mathtoolsset{showonlyrefs}

\begin{document}
\title[Algebraic curves $A^{\circ l}(x)-U(y)=0$ and 
arithmetic of orbits]{Algebraic curves $A^{\circ l}(x)-U(y)=0$  and 
arithmetic of orbits of rational functions
}
\author{F. Pakovich}
\thanks{
This research was partially supported by ISF Grant No. 1432/18}
\address{Department of Mathematics \\
Ben-Gurion University of the Negev \\
P.O.B. 653 Beer Sheva \\
8410501 Israel
}
\email{pakovich@math.bgu.ac.il}

\begin{abstract} We give a description of pairs of complex rational functions $A$ and $U$ of degree at least two such that for every $d\geq 1$ the 
algebraic curve \linebreak $A^{\circ d}(x)-U(y)=0$ has a factor of genus zero or one. In particular, we show that if  $A$ is not a ``generalized Latt\`es map'', then this 
condition is satisfied if and only if
there exists a rational function $V$ such that $U\circ V=A^{\circ l}$ for some $l\geq 1.$
We also prove a version of  the dynamical Mordell-Lang conjecture,
concerning intersections of
orbits of points from $\P^1(K)$ under iterates of $A$  
with the value set $U(\P^1(K))$, where $A$ and $U$ are rational functions defined over a number field $K.$

\end{abstract}


\maketitle

\begin{section}{Introduction}
In this paper we solve the following problem.

\bpr \l{pr1} Describe the pairs of complex rational functions $A$ and $U$  of degree at least two such that for every $d\geq 1$, the algebraic curve 
\be \l{1} A^{\circ d}(x)-U(y)=0\ee has an irreducible factor of genus zero or one.
\epr

The motivation for this problem comes from the  
arithmetic dynamics.
Spe\-ci\-fically, in \cite{jo}, the following problem was investigated: which rational functions $A$ defined over a number field $K$ have a $K$-orbit containing infinitely many distinct $m$th powers of elements from $K$ ? If such an orbit exists, then  for every $d\geq 1$, the algebraic curve 
\be \l{s}  A^{\circ d}(x)-y^m=0\ee has infinitely many $K$-points, implying by the Faltings theorem that it has a factor of genus zero or one. Thus, 
  a geometric counterpart of the initial arithmetic problem is to  describe rational functions $A$ such that all curves \eqref{s} have 
  a factor of genus zero or one. Now if instead of intersections of orbits of $A$ with the set of $m$th powers we consider intersections with the value set $U(\P^1(K))$ of an arbitrary rational function $U$, we arrive at  Problem \ref{pr1}.

The paper \cite{jo}, based on painstaking  calculations of the possible ramifications of rational functions $A$ such that every curve \eqref{s} has a factor of genus zero or one, provides a very 
explicit description of such functions. In contrast, our approach is 
more geometric and 
provides an answer in the general case in terms of semiconjugacies  and Galois coverings.
Notice that Problem \ref{pr1} is somewhat  similar to the following problem considered in the paper \cite{cur}: for which
rational functions $U$, there exists a sequence of rational functions $F_d,$ $d\geq 1,$ such that 
$\deg F_d\to \infty$, and for every $d\geq 1$,  
the curve
\be \l{seq} F_d(x)-U(y)=0\ee is irreducible and of genus zero. It was shown in \cite{cur} that $U$ satisfies this condition if and only if the Galois closure of the field extension $\C(z)/\C(U)$ has genus zero or one. 
Thus, this condition also holds for solutions of  Problem \ref{pr1} whenever curves \eqref{1} are irreducible. 
However, Problem \ref{pr1} is distinct from the problem considered in \cite{cur} in two important respects. First, curves \eqref{1} can be reducible. Secondly,  Problem \ref{pr1} asks for a description 
of all {\it pairs} $A$, $U$ such that curves \eqref{1} have an irreducible factor of genus zero or one, and not for a description of $U$ for which {\it some} $A$ with this property exists.

Let $A$ and $B$ be rational functions of degree at least two. Recall that 
the function $B$ is called semiconjugate to the function $A$
if there exists a non-constant rational function $X$
such that the diagram 
\be \l{ii1}
\begin{CD}
\C\P^1 @>B>> \C\P^1 \\
@VV X V @VV  X V\\ 
\f \C\P^1 @>A>> \f\C\P^1\
\end{CD}
\ee
commutes. 
 Semiconjugate rational functions  appear naturally in complex and arithmetic dynamics (see, e.g., the recent papers \cite{e}, \cite{ms}, \cite{pj}). They are also  
closely related to Problem \ref{pr1}. Indeed, 
since the commutativity of diagram \eqref{ii1} implies that 
$$ A^{\circ d}\circ X=X\circ B^{\circ d}, \ \ \  \ d\geq 1,$$  setting $U$ equal to $X$, we see that 
for every $d\geq 1$, curve \eqref{1} has a component of genus zero 
paramet\-ri\-zed by the rational functions $X$ and $B^{\circ d}.$ 
 
More generally, if $A,B$ and $X$ satisfy \eqref{ii1}, then  curves \eqref{1}   have a factor of genus zero for any rational function $U$, which is a ``compositional left factor''  
of the function $A^{\circ l}\circ X$ for some $l\geq 0,$ where by  a compositional left factor of a holomorphic map $f:\, R_1\rightarrow R_2$ between Riemann surfaces, we mean 
any holomorphic map  $g:\, R'\rightarrow R_2$ such that $f=g\circ h$ for some  holomorphic map  $h:\, R_1\rightarrow R'.$ 
Indeed, it follows from \eqref{ii1} and   $$A^{\circ l}\circ X=U\circ V$$ 
that 
$$A^{\circ l+k}\circ X=U\circ V\circ B^{\circ k}$$ for every $k\geq 0$, implying as above 
that the pair $A$, $U$ is a solution of Problem \ref{pr1}.  
In particular, setting $B=A$ and $X=z$ in \eqref{ii1}, we see that for every $d\geq 1$ curve \eqref{1} have a factor of genus zero whenever $U$  is
a  compositional left factor of some iterate  $A^{\circ l}$, $l\geq 1.$

Semiconjugate rational functions were studied at length in  the recent series of papers \cite{semi}, \cite{arn}, \cite{rec},  \cite{lattes}, \cite{fin} by using the theory of  orbifolds on Riemann surfaces. Our approach to Problem \ref{pr1} is based on the ideas and methods described in these papers. Roughly speaking, our main result  states  that, unless $A$ belongs to a special family of functions, all corresponding solutions $U$ of Problem \ref{pr1} 
can be obtained as described above from some {\it fixed} semiconjugacy \eqref{ii1}, where $X$ is a {\it Galois covering} that depends only on $A$.   
Moreover,  
 for ``most'' rational functions $A$, this Galois covering $X_A$ is equal simply to the identity map. In other words, a rational function $U$ is a solution of Problem \ref{pr1} 
 if and only if $U$ is 
a  compositional left factor of $A^{\circ l}$ for some $l\geq 1.$
 
To formulate our results explicitly we need several definitions.
Recall that an {\it orbifold} $\f O$ on $\C\P^1$  is a ramification function $\nu:\C\P^1\rightarrow \mathbb N$ which takes the value $\nu(z)=1$ except at a finite set of points. 
 We assume that considered orbifolds are {\it good}, meaning that we forbid $\f O$ to have exactly one point with $\nu(z)\neq 1$ or exactly two such points $z_1,$ $z_2$ with $\nu(z_1)\neq \nu(z_2).$ 
Let $f$ be a rational function and $\f O_1$, $\f O_2$  orbifolds with ramification functions $\nu_1$ and $\nu_2$. We say that  
$f:\,  \f O_1\rightarrow \f O_2$ is  {\it a covering map} 
between orbifolds
if for any $z\in \C\P^1$ the equality 
$$ \nu_{2}(f(z))=\nu_{1}(z)\deg_zf$$ holds.
In case   
the weaker condition 
$$ \nu_{2}(f(z))=\nu_{1}(z)\GCD(\deg_zf, \nu_{2}(f(z))$$ is satisfied,  
we say that $f:\,  \f O_1\rightarrow \f O_2$ is  {\it a  minimal holomorphic  map} 
between orbifolds.

In the above terms, {\it a  Latt\`es map} can be defined as a rational function $A$  such that $A:\f O\rightarrow \f O$ is a  covering self-map
for some orbifold $\f O$ (see \cite{mil2} for a classical  definition and for a proof of the equivalency
between two definitions).  
Following \cite{lattes}, we say 
 that a rational function $A$ is {\it a generalized Latt\`es map}  if there exists an orbifold $\f O$ 
 distinct from the non-ramified sphere such that  $A:\,  \f O\rightarrow \f O$ is  a  minimal holomorphic  map. 
Thus, $A$ is a  Latt\`es map if there exists an orbifold $\f O$ such that 
$$ \nu(A(z))=\nu(z)\deg_zA, \ \ \ \ z\in \C\P^1,$$
and $A$ is a  generalized Latt\`es map if there exists an orbifold $\f O$ such that 
\be \l{ee} \nu(A(z))=\nu(z)\GCD(\deg_zA,\nu (A(z))),  \ \ \ \ z\in \C\P^1.\ee 
Similar to ordinary  Latt\`es maps, generalized Latt\`es maps can be described in terms of semiconjugacies and group actions. In particular, the following statement is true (see \cite{lattes}): a rational function $A$ is a  generalized Latt\`es map if and only if 
there exist a compact Riemann surface $ {R}$ of genus zero or one,  a finite non-trivial group $\Gamma\subseteq Aut( {R})$, 
a group automorphism $\phi:\Gamma\rightarrow \Gamma,$ and  a holomorphic map $B:\, {R} \rightarrow   {R}$   
such that the diagram 
\be \l{xui+}
\begin{CD}
 {R} @>  B>>  {R} \\
@VV  \pi V @VV \pi  V\\ 
\C\P^1 @>A >> \ \ \C \P^1\, ,
\end{CD}
\ee
where $\pi:\, R\rightarrow R/\Gamma$ is the quotient map, commutes, 
and for any $\sigma\in \Gamma$ the equality
\be \l{homm2}  B\circ\sigma=\phi(\sigma)\circ B \ee holds.

We say that a rational function is {\it special} if  it is either a Latt\`es map or it is conjugate to $z^{\pm n}$ or $\pm T_n,$ where $T_n$ is the Chebyshev polynomial.  
For rational functions $A$ and $U$, we denote 
by $g_d=g_d(A,U)$, $d\geq 1,$  the minimal number $g$ such that curve \eqref{1} has a component of genus $g.$ In this notation our main result 
concerning Problem \ref{pr1} 
is 
following.

\bt \l{t1} Let  $A$ be a  non-special rational function of degree at least two.
 Then  there exist  a rational Galois covering $X_A$ and a rational function $B$ such that 
the diagram 
\be \l{udod}
\begin{CD}
\C\P^1 @>B>> \C\P^1 \\
@VV X_A V @VV  X_A V\\ 
\f \C\P^1 @>A>> \f\C\P^1\
\end{CD}
\ee
commutes, and for a rational function $U$ of degree at least two the sequence $g_d,$ $d\geq 1,$ is bounded 
 if and only if $U$ is a compositional left factor  of $A^{\circ l}\circ X_A $ for some $l\geq 0.$ 
Moreover, if $A$ is not  a generalized Latt\`es map, then  $g_d,$ $d\geq 1,$ is bounded   if and only if $U$ is a compositional left factor  of $A^{\circ l} $ for some $l\geq 1.$ 
\et 

Notice that our method provides an explicit description of 
the Galois covering $X_A$ appearing in Theorem \ref{t1} via the ``maximal'' orbifold $\f O$ for which \eqref{ee} is satisfied. In particular, the function $X_A$  
is defined by the function $A$ in a unique way up to natural isomorphisms.

Theorem \ref{t1}  can be illustrated as follows. A ``random'' rational function $A$ is  not a generalized Latt\`es map. Thus, a rational function $U$ is a solution of  Problem \ref{pr1} if and only if $U$ is a
compositional left factor of $A^{\circ l}$ for some $l\geq 1.$ 
A simple example of a generalized Latt\`es map 
is provided by any function of the form  $A=z^rR^n(z)$, where $R$ is a rational function, $n\geq 2$, $r\geq 1,$ and $\GCD(r,n)=1.$ Indeed, one can easily check that \eqref{ee} is satisfied for the orbifold 
$\f O$ defined by the conditions $$ \nu(0)=n, \ \ \ \nu(\infty)=n.$$ 
With a few exceptions, the rational function $A=z^rR^n(z)$ is not special, and diagram \eqref{udod} from Theorem \ref{t1} has the form 
$$
\begin{CD}
\C\P^1 @>z^rR(z^n)>> \C\P^1 \\
@VV z^n V @VV  z^n V\\ 
\f \C\P^1 @>z^rR^n(z)>> \f\C\P^1\, .
\end{CD}
$$
Thus, for such $A$ a rational function $U$ is a
solution of Problem \ref{pr1}   if and only if there exists $l\geq 0$ such that $U$ is a  
 compositional left factor of the function $(z^rR^n(z))^{\circ l}\circ z^n$.

Assume now that considered rational functions $A$ and $U$ are defined over a number field $K.$ 
As an application of our results, we prove a statement that,
concerning intersections of
orbits of points from $\P^1(K)$ under iterates of $A$  
with the value set $U(\P^1(K))$,
can be considered to be a version of  the dynamical Mordell-Lang conjecture. 
Recall that the dynamical Mordell-Lang conjecture 
states that if $f$ is an endomorphism of a quasiprojective variety $V$ over $\C$, then for any point $z_0\in V$ and any subvariety $W\subset V$ 
the set of indices $n$ such that the $X^{\circ n}(z_0)\in W$ is a finite union of arithmetic progressions (see \cite{bgt} and the bibliography therein). 
In particular, this implies that if the $f$-orbit of $z_0$ has an infinite intersection with a proper subvariety $W$, then its Zariski closure is contained in a finite union of  proper subvarieties, and therefore, it  cannot coincide with whole $V$.
Notice that singletons are considered as arithmetic progressions with the common difference equal zero, so 
any finite set is a union of arithmetic progressions.

A conjecture closely related to 
the dynamical Mor\-dell-Lang conjecture was proposed in \cite{ms} (see also \cite{abr}, \cite{swz}).
It states that 
if  $f$ is a dominant endomorphism of a quasiprojective variety $V$ defined over an algebraically closed field $K$ of characteristic zero
 for which there exists no non-constant rational function $g$ satisfying $g \circ  f = g$, then there is  a point $z_0 \in V(K)$ whose $f$-orbit is Zariski dense in $V$. 
This conjecture  is complementary to the Mor\-dell-Lang conjecture
in the sense that the former states that there is a point with the dense orbit, while the dynamical Mordell-Lang conjecture asserts that the orbit of such a point intersects any subvariety $W$ of $V$ 
in at most finitely many
points.

In this paper, we prove the following statement that is similar in spirit to the 
dynamical Mordell-Lang conjecture.

\bt \l{t2} 
Let $A$ and $U$ be rational functions of degree at least two defined over a number field $K$, and $x_0$ a point from $\P^1(K).$  
Then the set of indices $n$ such that $A^{\circ n}(x_0)\in U(\P^1(K))$ is a finite union of arithmetic progressions. 
Moreover, if  $A$ is not  a generalized Latt\`es map, then either the above set is finite, or $A^{\circ n}(x_0)$ belongs  to $U(\P^1(K))$ for all but finitely many $n$.  
\et 
The first part of Theorem \ref{t2} confirms a conjecture proposed  in the paper \cite{jo}\footnote{A 
proof of this conjecture  
was also announced by T. Hyde and M. Zieve in the ``Workshop on Interactions between Model Theory and Arithmetic Dynamics'' in July 2016 at the Fields Institute. To date, however, a complete proof of their results is still not available.}. On the other hand, the second part asserts that if $A$ is not  a generalized Latt\`es map, then a stronger conclusion holds.

The paper is organized as follows. In the second section, we present relevant definitions and some results concerning orbifolds, fiber products, and generalized Latt\`es maps, mostly proved in the papers \cite{semi}, \cite{lattes}. 
In the third section, 
using the lower bounds obtained in   \cite{cur} on the genus of algebraic curves of the form  $$F(x)-U(y)=0,$$ where $F$ and $U$ are rational functions,   we solve Problem \ref{pr1}. In particular, we prove Theorem \ref{t1}. We also  
solve Problem \ref{pr1} for special $A$. 
In fact, we consider a more general version of Problem \ref{pr1} in which $U$ is  allowed to be a holomorphic map 
$$U:R\rightarrow \C\P^1,$$ where $R$ is a compact Riemann surface, and instead of curves \eqref{1} the
fiber products of $U$ and $A^{\circ d},$ $d\geq 1$, are  considered.

Finally, in the fourth section we prove Theorem \ref{t2}, combining the results of the third section with some results 
for semiconjugate maps between algebraic curves which are interesting in their own right. 
We also provide an example illustrating constructions and results of this paper.

\end{section}

\begin{section}{Orbifolds and generalized Latt\`es maps
}

\begin{subsection}{Riemann surface orbifolds}

In this section, we recall basic definitions concerning orbifolds on Riemann surfaces (see \cite{mil}, Appendix E) and 
some results and constructions 
from the papers \cite{semi}, \cite{lattes}.
 We also prove some additional related results used later. 

{\it A Riemann surface  orbifold} is a pair $\f O=(R,\nu)$ consisting of a Riemann surface $R$ and a ramification function $\nu:R\rightarrow \mathbb N$ 
that takes the value $\nu(z)=1$ except at isolated points. 
For an orbifold $\f O=(R,\nu)$ 
 the {\it  Euler characteristic} of $\f O$ is the number
\be \l{euler} \chi(\f O)=\chi(R)+\sum_{z\in R}\left(\frac{1}{\nu(z)}-1\right),\ee
the set of {\it singular points} of $\f O$ is the set 
$$c(\f O)=\{z_1,z_2, \dots, z_s, \dots \}=\{z\in R \mid \nu(z)>1\},$$ and  the {\it signature} of $\f O$ is the set 
$$\nu(\f O)=\{\nu(z_1),\nu(z_2), \dots , \nu(z_s), \dots \}.$$ 
For   orbifolds $\f O_1=(R_1,\nu_1)$  and $\f O_2=(R_2,\nu_2)$, 
we write \be \l{elki} \f O_1\preceq \f O_2 \ee 
if $R_1=R_2$, and for any $z\in R_1$, the condition $\nu_1(z)\mid \nu_2(z)$ holds.
Clearly, \eqref{elki} implies that 
$$\chi(\f O_1)\geq \chi(\f O_2).$$

Let $\f O_1=(R_1,\nu_1)$  and $\f O_2=(R_2,\nu_2)$ be orbifolds
and let 
$f:\, R_1\rightarrow R_2$  be a holomorphic branched covering map. We say that $f:\,  \f O_1\rightarrow \f O_2$
is  a {\it covering map} 
{\it between orbifolds}
if for any $z\in R_1$, the equality 
\be \l{us} \nu_{2}(f(z))=\nu_{1}(z)\deg_zf\ee holds, where $\deg_zf$ is the local degree of $f$ at the point $z$.
If for any $z\in R_1$, 
the weaker condition 
\be \l{uuss} \nu_{2}(f(z))\mid \nu_{1}(z)\deg_zf\ee
is satisfied  instead of  \eqref{us},   we say that $f:\,  \f O_1\rightarrow \f O_2$ 
is a {\it holomorphic map}
 {\it between orbifolds}.

If $f:\,  \f O_1\rightarrow \f O_2$ is a covering map between orbifolds with compact $R_1$ and $R_2$, then  the Riemann-Hurwitz 
formula implies that 
\be \l{rhor} \chi(\f O_1)=d \chi(\f O_2), \ee
where $d=\deg f$. 
For holomorphic maps, the following statement is true (see \cite{semi}, Proposition 3.2). 

\bp \l{p1} Let $f:\, \f O_1\rightarrow \f O_2$ be a holomorphic map between orbifolds with compact $R_1$ and $R_2$.
Then 
\be \l{iioopp} \chi(\f O_1)\leq \chi(\f O_2)\,\deg f, \ee and the equality 
holds if and only if $f:\, \f O_1\rightarrow \f O_2$ is a covering map between orbifolds. \qed
\ep

Let $R_1$, $R_2$ be Riemann surfaces and 
$f:\, R_1\rightarrow R_2$ a holomorphic branched covering map. Assume that $R_2$ is provided with a ramification function $\nu_2$. To define a ramification function $\nu_1$ on $R_1$ so that $f$ would be a holomorphic map between orbifolds $\f O_1=(R_1,\nu_1)$ and $\f O_2=(R_2,\nu_2)$, 
we must satisfy condition \eqref{uuss}, and it is easy to see that
for any  $z\in R_1$, a minimum possible value for $\nu_1(z)$ is defined by 
the equality 
\be \l{rys} \nu_{2}(f(z))=\nu_{1}(z)\GCD(\deg_zf, \nu_{2}(f(z)).\ee 
In case \eqref{rys} is satisfied for  any $z\in R_1$, we 
say that $f:\,  \f O_1\rightarrow \f O_2$  is  a  {\it minimal holomorphic  map 
between orbifolds}. 
It follows from the definition that for any orbifold $\f O=(R,\nu)$ and holomorphic branched covering map $f:\, R^{\prime} \rightarrow R$, there exists a unique orbifold structure $\nu^{\prime}$ on $R^{\prime}$, such that 
$f$ becomes a minimal holomorphic map between orbifolds. We will denote the corresponding orbifold by $f^*\f O.$

Below we will use 
the following property of the association from $\f O$ to $f^*\f O$ 
(see \cite{semi}, 
Corollary 4.2).

\bp \l{indu2}  Let $f:\, R_1 \rightarrow R^{\prime}$ and $g:\, R^{\prime} \rightarrow R_2$ be holomorphic branched covering maps, and  $\f O_1=(R_1,\nu_1)$ and  $\f O_2=(R_2,\nu_2)$
orbifolds. Assume that \linebreak $g\circ f:\, \f O_1\rightarrow \f O_2$ is  a minimal holomorphic map (resp. a co\-vering map). Then  $f:\, \f O_1\rightarrow g^*\f O_2 $ and $g:\, g^*\f O_2\rightarrow \f O_2$ are minimal holomorphic maps (resp. covering maps). \qed
\ep 

A universal covering of an orbifold ${\f O}$
is a covering map between orbifolds \linebreak  $\theta_{\f O}:\,
\tt {\f O}\rightarrow \f O$ such that $\tt R$ is simply connected and $\tt {\f O}$ is non-ramified, that is, $\tt \nu(z)\equiv 1.$ 
If $\theta_{\f O}$ is such a map, then 
there exists a group $\Gamma_{\f O}$ of conformal automorphisms of $\tt R$ such that the equality 
$$\theta_{\f O}(z_1)=\theta_{\f O}(z_2)$$ holds for $z_1,z_2\in \tt R$ if and only if $z_1=\sigma(z_2)$ for some $\sigma\in \Gamma_{\f O}.$ 
A universal covering exists and 
is unique up to a conformal isomorphism of $\tt R$ whenever 
$\f O$ is {\it good}, that is,  distinct from the Riemann sphere with one ramified point or with two ramified points $z_1,$ $z_2$ such that $\nu(z_1)\neq \nu(z_2)$.   
 Furthermore, 
$\tt R$ is the unit disk $\mathbb D$ if and only if $\chi(\f O)<0,$ $\tt R$ is the complex plane $\C$ if and only if $\chi(\f O)=0,$ and $\tt R$ is the Riemann sphere $\C\P^1$ if and only if $\chi(\f O)>0$ (see e.g. \cite{fk}, Section IV.9.12).
Below we will always assume that considered orbifolds are good.
Abusing  notation, we will use the symbol $\tt {\f O}$ both for the
orbifold and for the  Riemann surface  $\tt R$.

Covering maps between orbifolds lift to isomorphisms between their univer\-sal coverings.
More generally, the following proposition is true (see \cite{semi}, Proposi\-tion 3.1).

\bp \l{poiu} Let $f:\,  \f O_1\rightarrow \f O_2$ be a holomorphic map between orbifolds. Then for any choice of $\theta_{\f O_1}$ and $\theta_{\f O_2}$, there exist 
a holomorphic map $F:\, \tt {\f O_1} \rightarrow \tt {\f O_2}$ and 
a homomorphism $\phi:\, \Gamma_{\f O_1}\rightarrow \Gamma_{\f O_2}$, such that the diagram 
\be \l{dia2}
\begin{CD}
\tt {\f O_1} @>F>> \tt {\f O_2}\\
@VV\theta_{\f O_1}V @VV\theta_{\f O_2}V\\ 
\f O_1 @>f >> \f O_2\ 
\end{CD}
\ee
is commutative, and 
for any $\sigma\in \Gamma_{\f O_1}$, the equality
\be \l{homm}  F\circ\sigma=\phi(\sigma)\circ F \ee holds.
The map $F$ is defined by $\theta_{\f O_1}$, $\theta_{\f O_2}$, and $f$  
uniquely up to a transformation 
$F\rightarrow g\circ F,$ where $g\in \Gamma_{\f O_2}$. 
In the other direction, for any holomorphic map $F:\, \tt {\f O_1} \rightarrow \tt {\f O_2}$  that satisfies \eqref{homm} for some homomorphism $\phi:\, \Gamma_{\f O_1}\rightarrow \Gamma_{\f O_2}$
there exists a uniquely defined  holomorphic map between orbifolds $f:\,  \f O_1\rightarrow \f O_2$ such that diagram \eqref{dia2} is commutative.
The holomorphic map $F$ is an isomorphism if and only if $f$ is a covering map between orbifolds. \qed
\ep

 With each holomorphic map $f:\, R_1\rightarrow R_2$ between compact Riemann surfaces, 
one can associate two orbifolds $\f O_1^f=(R_1,\nu_1^f)$ and 
$\f O_2^f=(R_2,\nu_2^f)$ in a natural way, setting $\nu_2^f(z)$  
equal to the least common multiple of local degrees of $f$ at the points 
of the preimage $f^{-1}\{z\}$, and $$\nu_1^1(z)=\frac{\nu_2^f(f(z))}{\deg_zf}.$$ By construction, 
 $$f:\, \f O_1^f\rightarrow \f O_2^f$$ 
is a covering map between orbifolds.
It is easy to  see that this covering map is minimal in the following sense. For any covering map   $f:\, \f O_1\rightarrow \f O_2$, we have:
$$ \f O_1^f\preceq \f O_1, \ \ \ \f O_2^f\preceq \f O_2.$$
The orbifolds $\f O_1^f$ and 
$\f O_2^f$ 
are  good (see \cite{semi}, Lemma 4.2).

\vskip 0.2cm

\bt \l{xriak+} Let $f:\, R_1\rightarrow R_2$ be a holomorphic map between compact Riemann surfaces and $\f O=(R_2,\nu)$ an orbifold. 
Then $f$ is a compositional left factor of  $\theta_{\f O}$ if and only if $\f O^f_2\preceq \f O$. Furthermore, for any 
decomposition  $\theta_{\f O}=f\circ \psi,$ where $\psi:\t {\f O}\rightarrow R_1$ is a holomorphic map, the equality $\psi=\theta_{f^*\f O}$ holds, and the map 
$f:\, f^*\f O\rightarrow \f O$ is a covering map between orbifolds. In particular, $\theta_{{\f O}_2^f}=f \circ  \theta_{{\f O}_1^f}.$ 
\et 
\pr Since $\f O_2^{\theta_{\f O}}=\f O$, the ``only if'' part follows from the chain rule. 
In the other direction, let $\psi$ be the analytic continuation of
$f^{-1}\circ \theta_{\f O}$, where $f^{-1}$ is a germ of the function inverse to $f$. It follows easily from the definitions and the condition $\f O^f_2\preceq \f O$ that 
$\psi$ has no ramification. Therefore, since $\t{\f O}$ is simply connected, $\psi$ is single-valued, and   $\theta_{\f O}=f\circ \psi$.

Finally, it follows from the equality $\theta_{\f O}=f\circ \psi$ by Proposition \ref{indu2} that
$$f:\, f^*\f O \rightarrow \f O, \ \ \ \ \psi:\, \t{\f O}\rightarrow f^*\f O$$ are  covering maps between orbifolds, implying that 
 $\psi=\theta_{f^*\f O}$, since $\t{\f O}$ is non-ramified and  simply-connected.
  In particular, if  $\f O=\f O^f_2,$ then  $f^*\f O^f_2=\f O_1^f,$ so that $\psi= \theta_{{\f O}_1^f}.$
\qed

\bc \l{newco} Let $C$ and $W$ be rational functions, and $\f O=(\C\P^1,\nu)$ an orbifold such that  
 $\f O^W_2\preceq \f O$. Then any compositional left factor $U$ of $C\circ W$ is a compositional left factor of $C\circ \theta_{\f O}.$
In particular, any compositional left factor of $C\circ W$ is a compositional left factor of  $C\circ \theta_{\f O_2^W}.$ 
\ec
\pr Indeed, the equalities $C\circ W=U\circ V$ and $\theta_{\f O}=W\circ \psi$ imply the 
equality $$C\circ \theta_{\f O}=U\circ (V\circ \psi). \eqno{\Box}$$

\bc \l{baran} Let $f:R\rightarrow \C\P^1$ be a holomorphic map between compact Riemann surfaces. Then 
$\chi(\f O_2^f)>0$ implies that $g(R)=0$. On the other hand, $\chi(\f O_2^f)= 0$ implies that  $g(R)\leq 1$. 
\ec 
\pr If  $\chi(\f O_2^f)>0$, then $\t{ {\f O}_2^f}=\C\P^1$. Thus, by Theorem \ref{xriak+},  
 $\theta_{{\f O}_1^f} : \C\P^1\rightarrow R$ is a holomorphic map,   implying  
that $g(R)=0$. Similarly, if $\chi(\f O_2^f)=0$, then 
 $\theta_{{\f O}_1^f} : \C \rightarrow R$ is a holomorphic map, implying that $g(R)\leq 1,$ since otherwise 
 lifting $\theta_{{\f O}_1^f}$ to a map between universal coverings (in the usual sense) would result in a contradiction with the Liouville theorem.\qed

\bc \l{xerez} Let $f:\, R\rightarrow \C\P^1$ be a holomorphic map between compact Riemann surfaces. Assume that $\f O_2^f$ is defined by the conditions
\be \l{los1} \nu_2^f(0)=n, \ \ \ \ \nu_2^f(\infty)=n.\ee Then $g(R)=0$, and $A=z^{n}\circ \mu$ 
for some M\"obius transformation $\mu$.  On the other hand, if $\f O_2^f$ is defined by the conditions
\be \l{los2} \nu_2^f(-1)=2, \ \ \ \ \nu_2^f(1)=2, \ \ \ \  \nu_2^f(\infty)=n, \ee  
then $g(R)=0$, and either $$f=\frac{1}{2}\left(z^n+\frac{1}{z^n}\right)\circ \mu,$$
or $f=\pm T_n\circ \mu$ for some M\"obius transformation $\mu.$
\ec
\pr Since by Theorem \ref{xriak+} the map $f$ is a compositional left factor of $\theta_{{\f O}_2^f}$, 
and the universal coverings for orbifolds given by \eqref{los1} and \eqref{los2} are rational functions 
\be \l{ffuu} Z_n=z^n, \ \ \ \ \ \ D_n=\frac{1}{2}\left(z^n+\frac{1}{z^n}\right)\ee correspondingly, the statement follows from the well-known fact that any 
compositional left factor of $Z_n$ has the form $Z_d\circ \mu$ for some M\"obius transformation $\mu$ and $d\vert n,$   
while any compositional  left factor of $D_n$ has the form $\pm T_d \circ \mu$ or $D_d\circ \mu$  
for some M\"obius transformation $\mu$ and $d\vert n$   (see e.g. \cite{gen}, Subsection 4.1 and 4.2). \qed

\end{subsection}

\begin{subsection}{Fiber products}
Let $f:\, C_1\rightarrow C$ and $g:\, C_2\rightarrow C$ be holomorphic maps between compact Riemann surfaces. 
The collection
\be \l{nota} (C_1,f)\times_C (C_2,g)=\bigcup\limits_{j=1}^{n(f,g)}\{R_j,p_j,q_j\},\ee 
where $R_j$ are compact Riemann surfaces provided with holomorphic maps
$$p_j:\, R_j\rightarrow C_1, \ \ \ q_j:\, R_j\rightarrow C_2, \ \ \ 1\leq j \leq n(f,g),$$
is called the {\it fiber product} of  $f$ and $g$ if \be \l{pes} f\circ p_j=g\circ q_j, \ \ \ 1\leq j \leq n(f,g),\ee 
and for any holomorphic maps $p:\, R\rightarrow C_1,$  $q:\, R\rightarrow C_2$
between compact Riemann surfaces satisfying 
$$f\circ p=g\circ q$$ there exist a uniquely defined  index $j$ and 
a holomorphic map $w:\, R\rightarrow R_j$ such that
 $$ p= p_j\circ  w, \ \ \ q= q_j\circ w.$$ 
The fiber product exists and is defined in a unique way up to natural isomorphisms.

In practical terms the fiber product is described by the following algebraic construction. Let us consider the algebraic curve 
\be \l{ccuurr} E=\{(x,y)\in C_1\times C_2 \, \vert \,  f(x)=g(y)\}.\ee
Let us denote by $V_j,$ $1\leq j \leq n(f,g)$,  irreducible components of $E$, by 
$R_j$, \linebreak  $1\leq j \leq n(f,g)$, their desingularizations, 
 and by $$\pi_j: R_j\rightarrow V_j, \ \ \ 1\leq j \leq n(f,g),$$ the desingularization maps.
Then the compositions  $$x\circ \pi_j:R_j\rightarrow C_1, \ \ \ y\circ \pi_j:R_j\rightarrow C_2, \ \ \ 1\leq j \leq n(f,g),$$ 
extend to holomorphic maps
$$p_j:\, R_j\rightarrow C_1, \ \ \ q_j:\, R_j\rightarrow C_2, \ \ \ 1\leq j \leq n(f,g),$$
and the collection $\bigcup\limits_{j=1}^{n(f,g)}\{R_j,p_j,q_j\}$ is the fiber product of $f$ and $g$.

Abusing notation we will call the Riemann  surfaces $R_j$,  $1\leq j \leq n(f,g),$ irreducible components of the fiber product of $f$ and $g$. 
The number of irreducible components  $n(f,g)$ satisfies the inequality 
\be \l{ineq} n(f,g)\leq \GCD(\deg f,\deg g).\ee 
Indeed, the degree of every map $$h_j=f\circ p_j=g\circ q_j, \ \ \ 1\leq j \leq n(f,g),$$ is divisible by  $\LCM(\deg f,\deg g)$. 
On the other hand, calculating the degrees of projections of  curve 
\eqref{ccuurr}, we see that  
\be \l{ii} \sum_{j}\deg p_j=  \deg g,\ \ \ \ \sum_{j}\deg q_j= \deg f,\ee
implying that $$\sum_{j=1}^{n(f,g)}h_j= \deg f \deg g.$$ 
Therefore, \eqref{ineq} holds.

\bt \l{sum} 
Let  $f:\, C_1\rightarrow C$, $g:\, C_2\rightarrow C$, and $u:\, C_3\rightarrow C_2$ be holo\-mor\-phic maps between compact Riemann surfaces. Assume that 
$$(C_1,f)\times_C (C_2,g)=\bigcup\limits_{j=1}^{n(f,g)}\{R_j,p_j,q_j\}$$
and 
$$(R_j, q_j)\times_{C_2}(C_3,u) =\bigcup\limits_{i=1}^{n(u,q_j)}\{R_{ij},p_{ij}, q_{ij}\},\ \ \ \ \ 1\leq j \leq n(f,g).$$ 
Then 
$$ (C_1,f)\times_C(C_3, g\circ u)=\bigcup\limits_{j=1}^{n(f,g)}\bigcup\limits_{i=1}^{n(u,q_j)}\{R_{ij},p_j\circ p_{ij},q_{ij}\}.$$ 
\et 
\pr 
It is clear that for $j,$ $1\leq j \leq n(f,g),$ and $i,$ $1\leq i \leq n(u,q_j)$, 
the diagram 
\be 
\begin{CD}
R_{ij} @>p_{ij}>> R_j  @>p_j>> C_1 \\
@VV q_{ij} V @VV q_j V @VV f V \\ 
C_3 @>u >> C_2  @>g >> C 
\end{CD}
\ee
commutes, so that  
 $$(g\circ u)\circ q_{ij}=f\circ (p_j\circ p_{ij}).$$

Assume now that $p$ and $q$ are holomorphic maps  between compact Riemann surfaces such that   
$$(g\circ u)\circ q=f\circ p.$$ By the universality property of the fiber product of $g$ and $f$, this equality  implies that 
$$u\circ q=q_j\circ w, \ \ \ \ p=p_j\circ w$$ for some  index $j$ and holomorphic map $w$. In turn, by the universality property of the fiber product of $u$ and $q_j$,   the first from these equalities implies  that 
$$q=q_{ij}\circ \t w, \ \ \ \ w=p_{ij}\circ \t w$$ for some   index $i$ and holomorphic map $\t w$.
Thus, $$p=p_j\circ p_{ij}\circ \t w, \ \ \ \ q=q_{ij}\circ \t w.\eqno{\Box}$$

\bc \l{sum+} In the above notation, the fiber products
$(C_1,f)\times_C (C_2,g)$ and $ (C_1,f)\times_C(C_3, g\circ u)$  have the same number of irreducible components if and only if for every $j,$ $1\leq j \leq n(f,g),$ the fiber product 
$(R_j, q_j)\times_{C_2}(C_3,u)$ has a unique irreducible component.  \qed 
\ec

\bc \l{volk} 
 Let  $R$ be a compact Riemann surface, $U:R\rightarrow \C\P^1$ a holomorphic map, and 
$A$ a rational function. Then there exists $d_0\geq 1$ such that \be \l{pk} n(A^{\circ d},U)=n(A^{\circ d_0},U)\ee 
for all $d\geq d_0$.

\ec
\pr 
Clearly, Theorem \ref{sum} implies that for every $d\geq 1$ the inequality $$n(A^{\circ (d+1)},U)\geq n(A^{\circ d},U)$$ holds. On the other hand, by \eqref{ineq}, for every $d\geq 1$  we have:
$$n(A^{\circ d},U)\leq \GCD(\deg A^{\circ d},\deg U)\leq \deg U.$$
Therefore, there exists $d_0\geq 1$ such that \eqref{pk} holds 
for all $d\geq d_0$. \qed

\end{subsection}

\begin{subsection}{Functional equations and orbifolds}
Orbifolds  $\f O_1^f$ and 
$\f O_2^f$ defined above are useful for the study of the functional equation
\be \l{m} f\circ p=g\circ q, \ee 
where  
$$p:\, R\rightarrow C_1,  \ \ \ \ f:\, C_1\rightarrow \C\P^1,\ \ \ \ q:\, R\rightarrow C_2, \ \ \ \ g:\, C_2\rightarrow \C\P^1$$ 
are holomorphic maps between compact Riemann surfaces.

We say 
that holomorphic maps $p:\, R\rightarrow C_1$ and $q:\, R\rightarrow C_2$ {\it have no non-trivial common compositional  right factor} if 
the equalities 
$$ p= \tt p\circ  w, \ \ \ q= \tt q\circ w,$$ where $w:\, R \rightarrow {\tt R}$, $\tt p:\, {\tt R}\rightarrow C_1$, $\tt q:\,  {\tt R}\rightarrow C_2$ are holomorphic maps between compact Riemann surfaces, imply that $\deg w=1.$ If such $p$ and $q$ satisfy \eqref{m}, then  by the universality property of the fiber product 
\be \l{upfp} (C_1,f)\times_{\C\P^1} (C_2,g)=\bigcup\limits_{j=1}^{n(f,g)}\{R_j,p_j,q_j\},\ee the equalities 
 $$p=p_j\circ w, \ \ \ \ q=q_j\circ w$$ hold for some $j$, $1\leq j \leq n(f,g),$ and an isomorphism $w:\, R_j\rightarrow R_j.$ 

A solution $f,p,g,q$ of \eqref{m} is called {\it good} if the fiber product of $f$ and $g$ has a unique component, 
and $p$ and $q$ have no non-trivial common compositional right factor. Thus, good solutions  correspond to fiber products \eqref{upfp} 
with $n(f,g)=1.$ 
In this notation, the following statement holds (see \cite{semi}, Theorem 4.2).

\bt \l{goodt} Let $f,p,g,q$ be a good solution of \eqref{m}.
Then the commutative diagram 
\be 
\begin{CD}
\f O_1^q @>p>> \f O_1^f\\
@VV q V @VV f V\\ 
\f O_2^q @>g >> \f O_2^f\ 
\end{CD}
\ee
consists of minimal holomorphic  maps between orbifolds. \qed  
\et 
 
Of course, vertical arrows in the above diagram are
holomorphic maps simply by definition. The meaning of the theorem is that the branching of $f$
and $q$ defines, to a certain extent, the branching of $g$ and $p$ and vice versa.

\vskip 0.2cm

Below we will use the following criterion (see \cite{semi}, Lemma 2.1).

 \bl \l{good} A 
solution $f,p,g,q$ of \e{m} is good whenever 
any two of the following three conditions are satisfied:

\begin{itemize}
\item the fiber product of $f$ and $g$ has a unique component,
\item $p$ and $q$ have no non-trivial common compositional right factor,
\item $ \deg f=\deg q, \ \ \ \deg g=\deg p.$  \qed
\end{itemize}
\el

\end{subsection}

\begin{subsection}{Generalized Latt\`es maps}\l{gelatt}

Most of orbifolds considered in this paper are defined on $\C\P^1$. 
For such orbifolds, we will omit  the Riemann surface $R$ in the definition of $\f O=(R,\nu)$,
meaning that $R=\C\P^1.$ 
Signatures of orbifolds on $\C\P^1$ with non-negative Euler characteristics and corresponding $\Gamma_{\f O}$ and $\theta_{\f O}$ can be described explicitly as follows. If $\f O$ is an orbifold distinct from the non-ramified sphere, then 
$\chi(\f O)=0$  if and only if the signature of $\f O$ 
belongs to the list
\be \l{list}\{2,2,2,2\} \ \ \ \{3,3,3\}, \ \ \  \{2,4,4\}, \ \ \  \{2,3,6\}, \ee and $\chi(\f O)>0$  if and only if
 the signature of $\f O$  belongs to the  list 
 \be \l{list2} \{n,n\}, \ \ n\geq 2,  \ \ \ \{2,2,n\}, \ \ n\geq 2,  \ \ \ \{2,3,3\}, \ \ \ \{2,3,4\}, \ \ \ \{2,3,5\}.\ee
Groups $\Gamma_{\f O}\subset Aut(\C)$ corresponding to orbifolds $\f O$ with signatures \eqref{list}  
are generated by translations of $\C$ by elements of some lattice $L\subset \C$ of rank two and the rotation $z\rightarrow  \v z,$ where $\v$ is an $n$th root of unity with $n$ equal to 2,3,4, or 6, such that  $\v L=L$ (see  \cite{mil2}, or \cite{fk}, 
Section IV.9.5).  
 Accordingly, the functions $\theta_{\f O}$ 
may be written in terms of the  corresponding
Weierstrass functions as $\wp(z),$ $\wp^{\prime }(z),$ $\wp^2(z),$  and $\wp^{\prime 2}(z).$  
Groups $\Gamma_{\f O}\subset Aut(\C\P^1)$ corresponding to   orbifolds $\f O$ with signatures \eqref{list2} are the well-known finite subgroups 
 $C_n,$  $D_{2n},$  $A_4,$ $S_4,$ $A_5$ of $Aut(\C\P^1)$, and the functions $\theta_{\f O}$ are Galois coverings of $\C\P^1$ by $\C\P^1$ of degrees 
$n$, $2n,$ $12,$ $24,$  $60,$ calculated for the first time by Klein in \cite{klein}.

A {\it Latt\`es map} can be defined as a rational function $A$ of degree at least two such that $A:\f O\rightarrow \f O$ is a  covering self-map
for some orbifold $\f O$ (see \cite{mil2}). Thus, $A$ is a Latt\`es map  if there exists an orbifold $\f O$
such that  for any $z\in \C\P^1$ the equality 
\be \l{uu0} \nu(A(z))=\nu(z)\deg_zA\ee holds. By  formula \eqref{rhor}, such $\f O$ necessarily satisfies $\chi(\f O)=0.$ Furthermore, 
for a given function $A$ there might be at most one  orbifold such that \eqref{uu0} holds (see \cite{mil2} and \cite{lattes}, Theorem 6.1).  

Following \cite{lattes}, we say that a rational function $A$ of degree at least two is 
a {\it generalized Latt\`es map} if there exists an orbifold $\f O$, 
 distinct from the non-ramified sphere, 
such that  $A:\f O\rightarrow \f O$ is a minimal holomorphic self-map between orbifolds; that is, for any $z\in \C\P^1$, the equality 
\be \l{uu} \nu(A(z))=\nu(z)\GCD(\deg_zA,\nu (A(z)))\ee holds.  By inequality \eqref{iioopp}, such $\f O$ satisfies $\chi(\f O)\geq 0$.
Since \eqref{uu0} implies \eqref{uu}, any ordinary Latt\`es map 
is a generalized Latt\`es map. 
Notice that if $\f O$ is the non-ramified sphere, then condition \eqref{uu} trivially holds for any rational function $A$. We say that a rational function is {\it special} if  it is either a Latt\`es map, or is conjugate to $z^{\pm n}$ or $\pm T_n,$ where $T_n$ is the Chebyshev polynomial.

In general, for a given function $A$ there might be several orbifolds $\f O$ satisfying \eqref{uu}, and even infinitely many such orbifolds. For example,  $z^{\pm d}:\f O\rightarrow \f O$ is a minimal holomorphic map for any $\f O$ defined by the conditions \be \l{raz} 
\nu(0)=\nu(\infty)=n, \ \ \ \ n\geq 2, \ \ \ \  \GCD(d,n)=1,\ee and  $\pm T_{d}:\f O\rightarrow \f O$ is a minimal holomorphic map  for any $\f O$ defined by  the conditions \be \l{dva} 
\nu(-1)=\nu(1)=2, \ \ \ \nu(\infty)=n,\ \ \ \ n\geq 1, \ \ \ \  \GCD(d,n)=1.\ee 
For odd $d$, additionally, $\pm T_{d}:\f O\rightarrow \f O$ is a minimal holomorphic map 
for $\f O$ defined 
by 
\be \l{tri} 
\nu(1)=2, \ \ \ \nu(\infty)=2,\ee 
or
\be \l{chet} 
\nu(-1)=2, \ \ \ \nu(\infty)=2.\ee   Nevertheless, the following statement holds (see \cite{lattes}, Theorem 1.2).

\bt \l{uni} 
Let $A$ be a rational function  of degree at least two not conjugate to $z^{\pm d}$ or $\pm T_d.$ Then there exists an orbifold $\f O_0^A$ such that $A:\, \f O_0^A\rightarrow \f O_0^A$
is a minimal holomorphic map between orbifolds, and for any orbifold $\f O$ such that 
$A:\, \f O\rightarrow \f O$ is a minimal holomorphic map between orbifolds, the relation $\f O\preceq \f O_0^A$ holds. Furthermore, $\f O_0^{A^{\circ l}}=\f O_0^A$ for any 
$l\geq 1$.  \qed
\et

Clearly, generalized Latt\`es maps are exactly rational functions for which the orbifold $\f O_0^A$ is distinct from the non-ramified sphere, completed
by the functions  $z^{\pm d}$ and $\pm T_d$ for which the orbifold $\f O_0^A$ is not defined. Furthermore,  ordinary Latt\`es maps are exactly  rational functions for which  $\chi(\f O_0^A)= 0$  (see \cite{lattes}, Lemma 6.4) and if $A$ is a  Latt\`es map, then the minimal holomorphic map $A:\f O_0^A\rightarrow \f O_0^A$ is a  covering map by Proposition \ref{p1}. Notice also that since a rational function $A$ is conjugate to $z^{\pm d}$ or $\pm T_d$ if and only 
 if some iterate $A^{\circ l}$ is conjugate to $z^{\pm ld}$ or $\pm T_{ld}$
 (see e.g. \cite{lattes}, Lemma 6.3),  Theorem \ref{uni} implies that $A$ is a generalized Latt\`es map if and only if some iterate $A^{\circ l}$  is a generalized Latt\`es map.

For exceptional functions $z^{\pm d}$ and $\pm T_d$, the 
orbifolds for which \eqref{uu} holds are described as follows
(see \cite{lattes}, Theorem 6.2).

\bt \l{214} Let $\f O$ be an orbifold distinct from the non-ramified sphere. 
\begin{enumerate}

\item The map  $z^{\pm d}:\f O\rightarrow \f O$, $d\geq 2,$ is a minimal holomorphic map between orbifolds if and only if $\f O$ is defined by conditions \eqref{raz}. 

\item The map $\pm T_{d}:\f O\rightarrow \f O$, $d\geq 2,$ is a minimal holomorphic map between orbifolds 
if and only if either $\f O$ is defined 
by conditions \eqref{dva},  or $d$ is odd and 
 $\f O$ is defined 
by conditions
\eqref{tri} or  
\eqref{chet}.  \qed

\end{enumerate}

\et

If $A$ is a generalized Latt\`es map, then $c(\f O_0^A)$ is a subset of the set  $c(\f O_2^A)$ consisting of critical values of $A$ unless $\deg A\leq 4.$ 
More generally, the following statement holds  (see \cite{lattes}, Lemma 6.6).

\bl \l{toch} Let $A$ be a rational function of degree at least five, and $\f O_1$, $\f O_2$ orbifolds distinct from the non-ramified sphere such that  $A:\f O_1\rightarrow \f  O_2$ is a minimal holomorphic map between orbifolds. Assume that $\chi(\f O_1)\geq 0$. Then    $c(\f O_2)\subseteq c(\f O_2^A)$. 
\el

\end{subsection}

\end{section}

\begin{section}{Algebraic curves $A^{\circ l}(x)-U(y)=0$ with components of low genus
}
In this section, we solve Problem \ref{pr1}. 
Our approach is based on the following theorem, which is a mild generalization of a result proved in \cite{cur}.

\bt \l{m2} 
Let $R$ be a compact Riemann surface and 
$W:R\rightarrow \C\P^1$ a holomorphic map of degree $n$. Then for any 
rational function $P$ of degree $m$
 such that the fiber product  of $P$ and $W$ consists of a unique component $ E$, the inequality  
\be \l{svko} \chi(E)\leq \chi(R)(n-1)-\frac{m}{42}\ee
holds, unless  $\chi(\f O_2^W)\geq 0.$ 
\et 
\pr 
The proof of Theorem \ref{m2} in the case $R=\C\P^1$  was given in \cite{cur}, Section 3. 
\linebreak 
The proof in the general case is obtained in the same way with appropriate modifications.
First of all, observe that if $q:\, E\rightarrow R$ is a holomorphic map  of degree $n$ bet\-ween compact Riemann surfaces, then 
\be  \l{cha} \chi(\f O_2^q)\geq  \chi(E)+\chi(R)(1-n).\ee
Indeed, it follows from \eqref{euler} that 
$$\chi(\f O_2^q)\geq  \chi(R)- c(q),$$ where $c(q)$ is the number of branch points of $q.$ On the other hand, since the number $c(q)$ is less than or equal to the number of points $z\in E$ where $\deg_z q>1$, the Riemann-Hurwitz formula 
$$\chi(E)=\chi(R)n-\sum_{z\in E}(\deg_zq-1) $$ implies  
that $$c(q)\leq \chi(R)n-\chi(E).$$ Therefore, \eqref{cha} holds.

Let $W\times_{\C\P^1} P=\{E,p,q\}.$ 
Since $$P:\, \f O_2^q\rightarrow \f O_2^W$$ is a minimal holomorphic map between orbifolds by Theorem \ref{goodt},  it follows 
from Proposition \ref{p1} that 
\be \l{burkl} \chi(\f O_2^q)\leq m \chi(\f O_2^W).\ee 
On the other hand,  \eqref{euler} implies that  if  $\chi(\f O)< 0,$ then in fact
\be \l{42} \chi(\f O)\leq -\frac{1}{42}\ee  
(where the equality is attained for the collection of ramification indices $(2,3,7)$). Therefore, if $\chi(\f O_2^W)<0$, then it follows 
from  \eqref{burkl} and \eqref{cha} that 
$$\chi(E)+\chi(R)(1-n)\leq -\frac{m}{42},$$
implying  \eqref{svko}. 
 \qed

Let us denote by $D= D\Big[R_d,A,W_d,h_d\Big]$
 an infinite commutative diagram 
\be \l{coma}
\dots \begin{CD}
@> \, >>  R_3 @> h_3 >> R_2 @> h_2 >> R_1 @>h_1 >> R_0\\
@.   @VV W_3 V @VV W_2 V @VV W_1 V @VV W_0 V \\
@> \, >>   \C\P^1 @>A >> \C\P^1\ @> A >> \C\P^1 @> A >> \C\P^1
\end{CD}
\ee 
consisting of holomorphic maps between compact Riemann surfaces. 
We say that $D$ is {\it good} if 
 for any $d_2>d_1\geq 0$,  the maps
\be \l{mapa} W_{d_1}, \ \ \ h_{d_1+1}\circ h_{d_1+2}\circ \dots \circ h_{d_2}, \ \ \ A^{\circ(d_2-d_1)}, \ \ \ W_{d_2}\ee form a good solution of equation \eqref{m}. Notice that if $D$ is good, then 
\be \l{step} \deg W_d=\deg W_0, \ \ \ \ d\geq 1,\ee by Lemma \ref{good}. 
We say that $D$ is 
 {\it preperiodic} if  there exist  $s_0\geq 0$ and $l\geq 1$ such 
that for any $d\geq s_0$ the Riemann surfaces $R_d$ and $R_{d+l}$ are isomorphic and 
\be \l{eqa} W_{d}=W_{d+l}\circ \alpha_{d}\ee for some isomorphism 
$$\alpha_{d}:R_d\rightarrow R_{d+l}.$$

Combined with the general properties of fiber products and generalized Latt\`es maps, Theorem \ref{m2} implies 
the following statement.

\bt \l{xyi} Let $D=D\Big[R_d,A,W_d,h_d\Big]$ be a diagram consisting of holomorphic maps of degree at least two.
Assume that $D$ is good and the sequence $g(R_d),$ $d\geq 0,$ is bounded. Then  $g(R_d)\leq 1$, $d\geq 0,$ and, unless $A$ is a Latt\`es map,
$g(R_d)=0$, $d\geq 0.$ Furthermore, $D$ is preperiodic and  
$A^{\circ l}:\f O_2^{W_{d}}\rightarrow \f O_2^{W_{d}}$ is a minimal holomorphic map between orbifolds for some $l\geq 1$ and all $d$ big enough. In particular, $A$ is a generalized Latt\`es map.
\et
\pr Since the sequence $g(R_d),$ $d\geq 0,$ is bounded from above, the sequence $\chi(R_d)$, $d\geq 0,$ is bounded from below. Therefore,  
applying Theorem \ref{m2} for $W=W_d$, $d\geq 0,$ and $P=A^{\circ j}$ with $j$ big enough, we conclude 
 that  \be \l{bur} \chi(\f O_2^{W_d})\geq 0,\ \ \ \ d\geq 0.\ee Hence,  $g(R_d)\leq 1$, $d\geq 0,$ 
by Corollary \ref{baran}. 
 
Let us show that the set of orbifolds $\f O_2^{W_{d}},$ $d\geq 0,$ contains only finitely many different orbifolds. Clearly, it is enough to show that the sequences $c(\f O_2^{W_{d}}),$ $d\geq 0,$ and $\nu(\f O_2^{W_{d}}),$ $d\geq 0,$ have only finitely many different elements.  
Since $D$ is good, it follows from The\-orem \ref{goodt} that
\be \l{cov} A^{\circ(d_2-d_1)}:\f O_2^{W_{d_2}}\rightarrow \f O_2^{W_{d_1}}\ee 
is a minimal holomorphic map between orbifolds  for any $d_2>d_1\geq 0$. In particular, 
 \be \l{cov+} A:\f O_2^{W_{d+1}}\rightarrow \f O_2^{W_{d}}, \ \ \ \ \ d\geq 0,\ee 
are minimal holomorphic maps. 
 Therefore,  if $\deg A>4$, then by Lemma \ref{toch}
every set $c(\f O_2^{W_{d}}),$ $d\geq 0,$ is a subset of the set $c(\f O_2^A)$, and hence, 
the sequence $c(\f O_2^{W_{d}}),$ $d\geq 0,$ has only finitely many different elements.
Moreover, this is true if  $\deg A\leq 4$. Indeed, the inequality $\deg A\geq 2$ implies 
the inequality $\deg A^{\circ 3}>4$, and hence, 
every set $c(\f O_2^{W_{d}}),$ $d\geq 0,$ is a subset of the set $c(A^{\circ 3})$, since 
$$A^{\circ 3}:\f O_2^{W_{d+3}}\rightarrow \f O_2^{W_{d}}, \ \ \ \ d\geq 0,$$ 
 also are  minimal holomorphic maps.
Finally, by \eqref{bur}, possible  signatures of the orbifolds $\f O_2^{W_{d}},$ $d\geq 0,$
belong  to lists \eqref{list} and \eqref{list2}, and it follows from equality \eqref{step} and 
Corollary \ref{xerez} that if  $\nu(\f O_2^{W_{d}})=\{n,n\}$, $n\geq 2,$ or  $\nu(\f O_2^{W_{d}})=\{2,2,n\}$, $n\geq 2,$ then  either $n=\deg W_0$ or $n=\deg W_0/2$. Therefore, the sequence $\nu(\f O_2^{W_{d}}),$ $d\geq 0,$ also has only finitely many different elements.

Since the set $\f O_2^{W_{d}},$ $d\geq 0,$ contains only finitely many different orbifolds,  
there exist an orbifold $\f O$ with $\chi(\f O)\geq 0$ and  a monotonically increasing  sequence $d_k\to \infty$ such that 
$\f O_2^{W_{d_k}}=\f O$, $k\geq 0.$ 
Moreover, by Theorem  \ref{xriak+}, the equalities 
\be \l{krot} \theta_{\f O}=W_{d_k}\circ \theta_{\f O_1^{W_{d_k}}}, \ \ \ k\geq 0,\ee
hold.
By the classification given in Subsection \ref{gelatt}, the group $\Gamma_{\f O}$  is finitely generated, 
and therefore, it has at most finitely many subgroups of any given index.
Since \eqref{krot} and \eqref{step} imply that 
every group $\Gamma_{\f O_1^{W_{d_k}}}$, $k\geq 0,$   has index $\deg W_0$
in $\Gamma_{\f O}$, 
we conclude that  the set of groups  $\Gamma_{\f O_1^{W_{d_k}}}$, $k\geq 0,$ contains only  finitely many different groups. Therefore, 
\be \l{pop} \Gamma_{\f O_1^{W_{d_{k_i}}}}=\Gamma_{\f O_1^{W_{d_{k_j}}}}\ee 
for some $k_i> k_j$ and hence 
$$\theta_{\f O_1^{W_{d_{k_i}}}}(x)= \alpha\circ\theta_{\f O_1^{W_{d_{k_j}}}}(y)$$ for some 
isomorphism  $\alpha:R_{d_{k_j}}\rightarrow R_{d_{k_i}},$
implying by \eqref{krot} the equality \be \l{is2} W_{d_{k_j}}=W_{d_{k_i}}\circ \alpha.\ee
Since $R_d$, $W_d$, $h_d$,  $d\geq 1$,   are defined 
by $R_{d-1}$ and $W_{d-1}$ in a unique way up to natural isomorphisms, \eqref{is2} implies that the preperiodicity 
condition holds for  $l=d_{k_i}-d_{k_j}$ and $s_0=d_{k_j}.$

Finally, setting $d_2=d+l$ and $d_1=d$ in \eqref{cov},  we see that if $d\geq s_0$, then 
$$A^{\circ l}:\f O_2^{W_d}\rightarrow \f O_2^{W_d}$$ is a minimal holomorphic map. Therefore, since the inequality $\deg W_d\geq
 2$ implies that $\f O_2^{W_d}$ cannot be the non-ramified sphere,
$A^{\circ l}$ is a
generalized Latt\`es map, and hence, $A$ is also a
generalized Latt\`es map. Moreover, unless $A$ is a Latt\`es map, $\chi(\f O_2^{W_d})>0$, $d\geq s_0.$ Therefore, $g(R_d)=0$, $d\geq s_0,$ by Corollary \ref{baran}, implying that $g(R_d)=0$ for all $d\geq 0,$
since $g(R_{d+1})\geq g(R_d),$ $d\geq 0.$
\qed

\vskip 0.2cm

Four theorems below provide a solution of Problem \ref{pr1}. 
The first theorem imposes no restrictions on the function $A$ and relates Problem \ref{pr1} with semiconjugacies.  
The other three provide more precise information for different classes of $A$. In particular, Theorem \ref{chain} implies Theorem \ref{t1} stated in the introduction.  
In fact, we examine a more general version of Problem \ref{pr1}, in which $U$ is allowed to be a holomorphic map 
$U:R\rightarrow \C\P^1,$ where $R$ is a compact Riemann surface, and instead of curves \eqref{1} the 
fiber products of  $U$ and $A^{\circ d},$ $d\geq 1$, are considered.

Let us denote 
by $g_d=g_d(A,U),$ $d\geq 1,$ the minimal number $g$ such that the fiber product of $U$ and $A^{\circ d}$  has a component of genus $g.$

\bt \l{chai} Let  $ R$ be a compact Riemann surface, $U: R\rightarrow \C\P^1$ a holomorphic map of degree at least two, and 
$A$ a rational function of degree at least two. 
 Then the sequence $g_d,$ $d\geq 1,$ is bounded 
 if and only if there exist a compact Riemann surface $ S$ of genus 0 or 1 and  holomorphic maps $F: S\rightarrow  S$ and 
$W: S\rightarrow \C\P^1$ such that the diagram 
\be \l{di}
\begin{CD}
 S @>F>> S \\
@VV W V @VV W V\\ 
\C\P^1 @>A^{\circ l_1} >> \ \ \C \P^1
\end{CD}
\ee
 commutes for some $l_1\geq 1,$  the fiber product of $W$
 and $A^{\circ l_1}$ consists of a unique component, 
$A^{\circ l_1}:\f O_2^{W}\rightarrow \f O_2^{W}$ is a minimal holomorphic map between orbifolds, 
and $U$ is a compositional left factor  of $A^{\circ l_2}\circ W $ for some $l_2\geq 0.$ 
In particular, if $A$ is not  a generalized Latt\`es map, then  $g_d,$ $d\geq 1,$ is bounded   if and only if $U$ is a compositional left factor  of $A^{\circ l} $ for some $l\geq 1.$ 
 
\et 
\pr To prove the sufficiency, observe that \eqref{di} 
and  \be \l{dii} A^{\circ l_2}\circ W=U\circ V\ee imply that 
$$A^{\circ (l_2+l_1k)}\circ W=U\circ V\circ F^{\circ k}, \ \ \  \ k\geq 0.$$  Therefore, for every $d\geq 1$, there exist holomorphic maps
of the form 
$$\phi_d=A^{\circ s_d}\circ W, \ \ \ \ \ \psi_d= V\circ F^{\circ r_d}, $$ where $s_d\geq 0,$ $r_d\geq 0,$
satisfying 
$$A^{\circ d}\circ \phi_d=U\circ \psi_d.$$ 
By the universality property of the fiber product, this implies that   
for every $d\geq 1$, there exist a component $\{E,p,q\}$ of  $A^{\circ d}\times U$ and a holomorphic map $w:S\rightarrow E$ such that 
$$\phi_d=p\circ w, \ \ \ \psi_d=q\circ w.$$ Clearly, for such $E$ we have: $$g(E)\leq g(S)\leq 1.$$

Let us now prove the necessity. 
Let $d_0$ be the number such that \eqref{pk} holds for all $d\geq d_0$, and 
let 
\be \l{deco} (\C\P^1,A^{\circ (d_0+k)})\times_{\C\P^1} (R,U)=\bigcup\limits_{j=1}^{s}\{R_{j,k},W_{j,k},H_{j,k}\}, \ \ \ \ \  \ k\geq 0,\ee
where $s=n(A^{\circ d_0},U).$
It follows from  the universality property of the fiber product and equality \eqref{pk} that  for every $k\geq 0$ and $j,$ $1\leq j \leq s$,  there exists a uniquely defined $j'$ such that 
$$H_{j',k+1}= H_{j,k}\circ h$$ for some holomorphic map $h:R_{j',k+1}\rightarrow R_{j,k}$, and without loss of generality we may assume that the numeration in \eqref{deco} is chosen in such a way that $j=j'.$ Thus, we can assume that for every $j,$  $1\leq j \leq s$, 
there exist holomorphic maps $h_{j,k}$, $k\geq 1$,  such that
$$H_{j,k}=H_{j,0}\circ h_{j,1}\circ h_{j,2}\circ \dots \circ h_{j,k}$$ and the diagram
\be \l{coma2}
\dots \begin{CD}
{{R}}_{j,3}@> h_{j,3} >> {{R}}_{j,2} @> h_{j,2} >> {{R}}_{j,1} @>h_{j,1} >> {{R}}_{j,0}\\
@VV W_{j,3} V @VV W_{j,2} V @VV W_{j,1} V @VV W_{j,0} V\\
 \C\P^1 @>A >> \C\P^1  @> A  >> \C\P^1 @> A   >> \C\P^1
\end{CD}
\ee commutes. Moreover, this diagram is good by Corollary  \ref{sum+}.
Finally, since 
\be \l{obvi} g({{R}}_{j,k+1})\geq g({{R}}_{j,k}), \ \ \ \ k\geq 0,\ee 
it follows from the boundness of  the sequence  $g_d,$ $d\geq 1,$ that  for at least one  $j,$ $1\leq j \leq s$,  the sequence 
 $g({{R}}_{j,k}),$ $k\geq 0,$  is  
 bounded. In particular,  for such $j$
 we can apply Theorem \ref{xyi} to diagram \eqref{coma2}, unless $\deg W_{j,0}= 1.$

By construction, for each $j$, $1\leq j \leq s,$ the diagram 
\be \l{coma3}
\dots \begin{CD}
{{R}}_{j,3}@> h_{j,3} >> {{R}}_{j,2} @> h_{j,2} >> {{R}}_{j,1} @>h_{j,1} >> {{R}}_{j,0}@> H_{j,0} >> R\\
@VV W_{j,3} V @VV W_{j,2} V @VV W_{j,1} V @VV W_{j,0} V @VV U V\\
 \C\P^1 @>A >> \C\P^1  @> A  >> \C\P^1 @> A   >> \C\P^1 @> A^{\circ d_0} >> \C\P^1\,
\end{CD}
\ee
commutes. Fix now  $j$ such that the sequence $g({{R}}_{j,k}),$ $k\geq 0,$ is bounded.  
 If $\deg W_{j,0}=1,$ then $R=\C\P^1$ and the equality 
$$ A^{\circ d_0}=U\circ H_{j,0}\circ W_{j,0}^{-1}$$ 
implies that $U$ is a compositional left factor  of $A^{\circ d_0} $. Therefore, in this case the theorem is true for $$S=\C\P^1, \ \ \  W=z, \ \ \ F=A,\ \ \ l_1=1,\ \ \ l_2=d_0.$$ 
On the other hand, if $\deg W_{j,0}\geq 2,$ then by Theorem \ref{xyi} there exist $l\geq 1$, $s_0\geq 0$, and 
 an isomorphism $$ \alpha :R_{j,s_0}\rightarrow R_{j,s_0+l}$$ such that  
$$ W_{j,s_0}=W_{j,s_0+l}\circ \alpha,$$ and $$A^{\circ l}:\f O_2^{W_{j,s_0}}\rightarrow  \f O_2^{W_{j,s_0}}$$ is a minimal holomorphic map. 
Thus, \eqref{di} holds for 
$$S=R_{j,s_0}, \ \ \ W=W_{j,s_0}, \ \ \ F= h_{j,s_0+1}\circ h_{j,s_0+2}\circ  \dots \circ  h_{j,s_0+l} \circ\alpha,$$ and $A$ is  a generalized Latt\`es map.
Finally, $U$ is a compositional left factor  of $A^{\circ l_2}\circ W $ for $l_2=d_0+s_0,$ 
since 
$$A^{\circ (d_0+s_0)}\circ W_{j,s_0}=U\circ H_{j,0}\circ h_{j,1}\circ h_{j,2}\circ  \dots \circ  h_{j,s_0}. \eqno{\Box} $$

Notice that in the proof of the sufficiency we did not use the assumptions that the fiber product of $W$ and $A^{\circ l_1}$ has one  component and $A^{\circ l_1}:\f O_2^{W}\rightarrow \f O_2^{W}$ is a minimal holomorphic map between orbifolds.
Thus, the theorem implies that if $U$ satisfy  \eqref{di} and \eqref{dii} for some $W,F$, and $V$, then it satisfies  \eqref{di} and \eqref{dii} for $W,F$, and $V$ that obey these conditions (cf. \cite{lattes}, Section 3).

\bt \l{chain} Let  $R$ be a compact Riemann surface, $U:R\rightarrow \C\P^1$ a holomorphic map of degree at least two, and 
$A$ a non-special rational function of degree at least two. 
 Then the sequence $g_d,$ $d\geq 1,$ is bounded 
 if and only if $R=\C\P^1$ and $U$ is a compositional left factor  of $A^{\circ l}\circ \theta_{\f O_0^A}$ for some $l\geq 1.$ 
In particular, if $A$ is not  a generalized Latt\`es map, then  $g_d,$ $d\geq 1,$ is bounded   if and only if $U$ is a compositional left factor  of $A^{\circ l} $ for some $l\geq 1.$ 
\et 
\pr If the sequence $g_d,$ $d\geq 1,$ is bounded, then by 
Theorem \ref{chai} there exist  a compact Riemann surface $ S$ of genus 0 or 1 and  a holomorphic map $W: S\rightarrow  \C\P^1$ such that $U$ is a compositional left factor of  $A^{\circ l_2}\circ W$ for some $l_2$, and 
$$A^{\circ l_1}:\f O_2^{W}\rightarrow \f O_2^{W}$$ is a minimal holomorphic map for some $l_1.$
On the other hand, since $A$ is not conjugate to $z^{\pm n}$ or $\pm T_n,$ the  orbifold $\f O_0^A$ is well-defined  
and 
\be \l{sh0} \f O_2^{W}\preceq \f O_0^{A^{\circ l_1}}=\f O_0^A,\ee  
 by
Theorem \ref{uni}. Thus, $U$ is  a compositional left factor of the holomorphic map $A^{\circ l_2}\circ  \theta_{\f O_0^A}$ by Corollary \ref{newco}. 
Moreover, since $A$ is not a Latt\`es map, $\chi(\f O_A^0)>0$ and $\theta_{\f O_0^A}$ and 
$U$ are rational functions.

In the other direction, since $\chi(\f O_0^A)>0$, Proposition \ref{poiu} implies that 
there exists a rational function $F$ such that the diagram 
\be \l{dia3}
\begin{CD}
\C\P^1 @>F>> \C\P^1\\
@VV\theta_{\f O_0^A}V @VV\theta_{\f O_0^A}V\\ 
\C\P^1 @>A >>\C\P^1
\end{CD}
\ee
commutes. Arguing now as in Theorem \ref{chai}, we 
conclude that if $U$ is a compositional left factor  of $A^{\circ l}\circ \theta_{\f O_0^A}$, then the sequence $g_d,$ $d\geq 1,$ is bounded.
 \qed

\bt  Let  $R$ be a compact Riemann surface, $U:R\rightarrow \C\P^1$ a holomorphic map of degree at least two, and 
$A$ a Latt\`es map. Then the sequence $g_d,$ $d\geq 1,$ is bounded 
 if and only if $U$ is a compositional left factor  of  
 $\theta_{\f O_0^A}.$ 
\et
\pr 
Arguing as in Theorem \ref{chain}, we conclude that if the sequence $g_d,$ $d\geq 1,$ is bounded, then $U$ 
 is a compositional left factor  of $A^{\circ l}\circ \theta_{\f O_0^A}$ for some $l\geq 1.$ 
Thus, to prove the necessity,  we must only show that if $A$ is a Latt\`es map, then any
compositional left factor of $A^{\circ l}\circ  \theta_{\f O_0^A}$, $l\geq 1,$ is a compositional left factor  of  
 $\theta_{\f O_0^A}.$ 
Recall that for a Latt\`es map $A$ the equality $\chi(\f O_0^A)=0$ holds and 
$A:\f O_0^A\rightarrow \f O_0^A$ is a covering map between orbifolds (see the remarks after Theorem \ref{uni}). 
Therefore, by Proposition \ref{poiu}, the function $F$ in diagram \eqref{dia2} is an isomorphism, implying that \eqref{dia2} takes the form 
\be \l{dia5}
\begin{CD}
\C @>F=az+b>> \C \\
@VV\theta_{\f O_0^A}V @VV\theta_{\f O_0^A}V\\ 
\C\P^1 @>A >>\C\P^1\,,
\end{CD}
\ee
where $a,b\in \C,$ $a\neq 0.$ Thus, for every $d\geq 1$ 
the equality 
\be \l{bbaa} \theta_{\f O_0^A}=A^{\circ d}\circ  \theta_{\f O_0^A}\circ (F^{-1})^{\circ d}\ee
holds, implying the necessary statement.

Assume now that $\theta_{\f O_0^A}=U\circ \psi$, where $\psi:\C \rightarrow R$ and $U:R\rightarrow \C\P^1$ are
holomorphic maps between Riemann surfaces. Since diagram \eqref{dia5} commutes, for every $d\geq 1$ the equality 
$$A^{\circ d}\circ \theta_{\f O_0^A}=U\circ (\psi \circ F^{\circ d})$$
holds, implying that the map $\psi_d:\C\rightarrow \C\P^1 \times R$ given by 
$$\psi_d:\,z\rightarrow (\theta_{\f O_0^A},\psi \circ F^{\circ d})$$ is a meromorphic 
parametrization of an irreducible component of the algebraic curve 
$$ E=\{(x,y)\in \C\P^1\times R \, \vert \,  A^{\circ d}(x)=U(y)\}.$$
Since an algebraic curve possessing a parametrization  by meromorphic functions on $\C$ has genus at most one, this implies 
that the sequence $g_d,$ $d\geq 1,$ is bounded. \qed

\bt  Let  $R$ be a compact Riemann surface, $U:R\rightarrow \C\P^1$ a holomorphic map of degree at least two, and 
$A$ a  rational function of degree at least two.
\begin{enumerate}

 \item If $A=z^{m},$  then the sequence $g_d,$ $d\geq 1,$ is bounded 
 if and only if $R=\C\P^1$ and $U=z^{s}\circ \mu,$ $s\geq 2,$ where $\mu$ is a M\"obius transformation, 
 
 \item If $A=T_m,$  then the sequence $g_d,$ $d\geq 1,$ is bounded 
 if and only if $R=\C\P^1$ and either $U=\pm T_s\circ \mu,$ $s\geq 2,$ or $$U=\frac{1}{2}\left(z^s+\frac{1}{z^s}\right)\circ \mu,\ \ \ \ s\geq 2,$$
where $\mu$ is a M\"obius transformation. 
 
\end{enumerate}
\et
\pr 
Let us consider the case $A=T_m$.  
In the case $A=z^m$ the proof is similar. For brevity, we will use the notation $D_n$ introduced in \eqref{ffuu}. Let us prove the necessity.
Applying Theorem \ref{chai} and keeping its notation, we observe first that if $\deg W=1$, then $U$ is a compositional left factor 
 of $T_{m^{l_2}}$. Therefore, since any compositional factor of $T_l$ has the form $T_d\circ \mu$ for some $d\vert l$ and M\"obius transformation $\mu$, in this case the statement is true.

Let us assume now that $\deg W>1$. 
Since $$T_{m^{l_1}}:\f O_2^{W}\rightarrow \f O_2^{W}$$ 
is a minimal holomorphic map between orbifolds, it follows from Theorem \ref{214} that  $\f O_2^W$ is defined by one of  conditions \eqref{dva}, \eqref{tri}, \eqref{chet}. 
Further, any compositional left factor  of the map $A^{\circ l_2}\circ W$ is a compositional left factor  of the map
$A^{\circ l_2}\circ  \theta_{\f O_2^W}$,  by Corollary \ref{newco}.
Since  
the universal coverings of the orbifolds given by \eqref{dva}, \eqref{tri}, \eqref{chet} are the 
function $D_n$ and the functions $-T_2$, $T_2$, correspondingly, this implies that 
$U$ is a
compositional left factor  either of the function 
$$T_{m^{l_2}}\circ D_n=D_{nm^{l_2}}$$ or of the function $$T_{m^{l_2}}\circ \pm T_2=\pm T_{2m^{l_2}}.$$
Since any 
compositional left factor of  $D_l$ has the form $\pm T_s \circ \mu$ or $D_s\circ \mu$  
for some M\"obius transformation $\mu$ and $s\vert l$, this proves the necessity.

Finally, since 
$$T_m\circ D_s=D_s\circ z^m$$ 
for any $s\geq 2,$ and $\pm T_s \circ \mu$ and $D_s\circ \mu$  are compositional left factors of  $D_s$, the sufficiency can be proved as in Theorem \ref{chai}. \qed

\end{section}

\begin{section}{Arithmetic of orbits of rational functions}
\begin{subsection}{Normalizations and definition fields}
Recall that for a non-constant holomorphic map between compact Riemann surfaces $X:C\rightarrow \C\P^1$,  
its  normalization $N_X$ is defined as a holomorphic map of the lowest possible degree
between compact Riemann surfaces  $N_X:\,S_X\rightarrow \C\P^1$,  such that $N_X$ is a Galois covering and
 \be \l{gops} N_X=X\circ H\ee for some  holomorphic map $H:\, S_X\rightarrow \C\P^1$. From the algebraic point of view the passage from $X$ to $N_X$ 
corresponds to the passage from the field extension $\f M(C)/X^*\C(z)$ to its Galois closure. 
The corresponding Galois group $G_X$ may be identified with the monodromy group of the covering  $X:C\rightarrow \C\P^1$, or,  in terms of  the normalization, with the group $Aut(S_X).$ 
The surface $S_X$ is defined up to isomorphism. For fixed $S_X$, the function $N_X$ is defined in a unique way while the function 
$H$ in \eqref{gops} 
is defined up to the change $H\rightarrow H\circ \mu,$ where $\mu\in Aut(S_X).$ 

\bt \l{itex} Let $C$ be a compact Riemann surface, and  $A:\C\P^1\rightarrow \C\P^1$, \linebreak $B:C\rightarrow C$, $X:C\rightarrow \C\P^1$ non-constant 
 holomorphic maps such that the diagram 
\be \l{oma}
\begin{CD}
C @> B >> C\\
@VV X V @VV X V  \\
 \C\P^1 @>A >> \C\P^1  
\end{CD}
\ee 
commutes and the fiber product of $A$ and $X$ consists of a unique component.
Then there exist holomorphic maps $F:S_X\rightarrow S_X$, $H:S_X\rightarrow C$, and a group automorphism $\phi:Aut(S_X)\rightarrow Aut(S_X)$ such that the equality $N_X=X\circ H$ holds,  the diagram 
\be 
\begin{CD} \l{dura}
S_X @> F >> S_X \\
@VV H V @VV H V\\ 
C @>B>> C \\
@VV X V @VV X V\\ 
\C\P^1 @>A >> \C\P^1\, ,  
\end{CD}
\ee
commutes,
and for any $\sigma\in Aut(S_X)$ the equality
\be   F\circ\sigma=\phi(\sigma)\circ F \ee holds.
\et 
\pr The proof is based on the following geometric description of $G_X$ and $S_X$ (see, e.g., \cite{f}, $\S$I.G or \cite{sot}, Section 2.2).
 Let $L$ be the $n$-fold fiber product of $X:C\rightarrow \C\P^1$ with itself, that is, the algebraic curve in $C^n$
 defined by the equation \be \l{uri} X(z_1)=X(z_2)=\dots =X(z_n).\ee
Let us denote by $\Delta$ 
the big diagonal of $C^n$, which consists of points where at least two coordinates  coincide,
and by $L_0$ the algebraic closure of $L\setminus \Delta$ in $L$. 
Then all irreducible components $V_1,$ $V_2, \dots ,V_r$ of $L_0$ are isomorphic, and the group $G_X$ can be identified with the subgroup 
of $S_n$ consisting of all permutations $\sigma\in S_n$, such that $$(x_{\sigma(1)}, x_{\sigma(2)}, \dots , x_{\sigma(n)})\in V_j$$ for some $j,$ $1\leq j \leq r,$ if and only if $$(x_1,x_2,\dots, x_n)\in V_j.$$ Furthermore, if 
$V$ is any irreducible component of $L_0$, and $\t{V}\xrightarrow{\eta} V$ is the desingularization map, then ${\t V}=S_X$ and 
the map $N_X$ is given by the composition
\be \l{pro} \t{V}\xrightarrow{\eta} {V}\xrightarrow{\pi_i}C\xrightarrow{X}\C\P^1,\ee
where  $\pi_i$ is the projection to any coordinate.

Define the maps $\f B:C^n\rightarrow C^n$, $\f A:(\C\P^1)^n\rightarrow (\C\P^1)^n$, and $\f X:C^n\rightarrow (\C\P^1)^n$ by the formulas 
$$\f A:\, (z_1,z_2,\dots, z_n)\rightarrow (A(z_1),A(z_2),\dots ,A(z_n)),$$ 
$$\f B:\, (z_1,z_2,\dots, z_n)\rightarrow (B(z_1),B(z_2),\dots ,B(z_n)),$$
$$\f X:\, (z_1,z_2,\dots, z_n)\rightarrow (X(z_1),X(z_2),\dots ,X(z_n)).$$
Clearly,  the diagram 
\be \l{kaba}
\begin{CD}
C^n @>\f B>> C^n \\
@VV\f X V @VV\f X V\\ 
(\C\P^1)^n @>\f A >> (\C\P^1)^n
\end{CD}
\ee
commutes, and by construction,  $L=\f X^{-1}(\Delta_0)$, where $\Delta_0$ is the usual diagonal in $(\C\P^1)^n$ which  consists of points where all coordinates  coincide.
 Therefore, since \linebreak $\f A(\Delta_0)=\Delta_0$, it follows from \eqref{kaba} that  $\f B(L)\subseteq L$. 

Let us show that \be \l{mor} \f B(L_0)\subseteq L_0.\ee 
Since the fiber product of $A$ and $X$ consists of a unique component, it follows from Lemma \ref{good} that the maps $X$ and $B$ have no non-trivial common compositional right factor, implying that 
$$[X^*\C(z),B^*\f M(C)]=\f M(C).$$
By the primitive element theorem, $B^*\f M(C)=\C[h]$ for some $h\in B^*\f M(C),$ so that 
$$\f M(C)=X^*\C(z)[h].$$ Clearly, this equality implies that for all but finitely many points  
$z\in \C\P^1$, the function $h$ takes $n$ distinct values on the set $X^{-1}\{z_0\}$.  
Since $h\in B^*\f M(C),$ this implies, in turn, that for all but finitely many points of 
$z\in \C\P^1$ the map $B$ takes $n$ distinct values on the set $X^{-1}\{z_0\}$, or equivalently 
that \eqref{mor} holds.

Let $V$ be an irreducible component of $L_0$. Then \eqref{mor} implies that $\f B(V)\subseteq V',$ where
$V'$ is another irreducible component of $L_0.$ Clearly, we can complete the diagram 
\be 
\begin{CD}
V @> \f B >> V' \\
@VV \pi_i V @VV \pi_i V\\ 
C @>B >> C \\
@VV X V @VV X V\\ 
\C\P^1 @>A >> \C\P^1\, ,  
\end{CD}
\ee
where  $\pi_i$ is the projection to any coordinate, 
to the diagram 
\be 
\begin{CD}
S_X @> F_0 >> S_X \\
@VV \eta V @VV \eta' V\\ 
V @> \f B >> V' \\
@VV \pi_i V @VV \pi_i V\\ 
C @>B >> C \\
@VV X V @VV X V\\ 
\C\P^1 @>A >> \C\P^1,  
\end{CD}
\ee
where $F_0:S_X\rightarrow S_X$ is a holomorphic map and
\be \l{nx} N_X=X\circ \pi_i\circ \eta'=X\circ \pi_i\circ \eta.\ee 
Since \eqref{nx} implies that  
$$\pi_i\circ \eta'=\pi_i\circ \eta\circ \mu$$
for some $\mu\in Aut(S_X)$, 
we conclude that diagram
\eqref{dura} commutes for $F=\mu\circ F_0$ and $H=\pi_i\circ \eta.$

Since $X\circ H$ is a Galois covering, it is easy to see that for any holomorphic map $F':S_X\rightarrow S_X$ that, along with $F$, satisfies \eqref{dura}, there exists $g\in Aut(S_X)$ such that  $F'=g\circ F$. 
In particular, for any $\sigma\in Aut(S_X)$ the equality 
$$ F\circ\sigma=g_{\sigma}\circ F$$ holds for some $g_{\sigma}\in Aut(S_X)$,  and it is easy to see 
that the map $$\phi: \sigma\rightarrow g_{\sigma}$$ is a group homomorphism. Finally, since $Aut(S_X)=G_X$,  
if $\Ker \phi\neq e,$    then there exists a non-identical permutation $\sigma\in G_X$ such that 
$$\f B(z_{\sigma(1)},z_{\sigma(2)},\dots,z_{\sigma(n)})=\f B(z_1,z_2,\dots, z_n)$$ for all $(z_1,z_2,\dots,z_n)\in V$.
Since this contradicts the fact that 
for all but finitely many points 
$z\in \C\P^1$, the map $B$ takes $n$ distinct values on the set $X^{-1}\{z_0\}$, we conclude that 
$\phi$ is an automorphism. \qed

Notice that since the only compact Riemann surfaces admitting self maps of degree $d>1$ are the Riemann sphere and tori, 
Theorem \ref{ite} implies that if $A,B$ and $X$ are rational functions of degree at least two such that diagram \eqref{oma} commutes and the fiber product of $A$ and $X$ consists of a unique component, then   $g(S_X)\leq 1$ 
(cf. \cite{arn}).

\bc \l{ite} Let $C$ be a compact Riemann surface and $A:\C\P^1\rightarrow \C\P^1$, \linebreak $B_1:C\rightarrow C$, $B_2:C\rightarrow C$,   $X:C\rightarrow \C\P^1$ non-constant 
holomorphic maps such that the diagrams
\be 
\begin{CD}
C @>B_1>> C \\
@VV X V @VV  X V\\ 
\C\P^1 @>A>> \C\P^1\,,
\end{CD}\ \ \ \ \\ \ \ \ \ \ \ 
\begin{CD}
C @>B_2>> C \\
@VV X V @VV  X V\\ 
\C\P^1 @>A>> \C\P^1
\end{CD}
\ee
 commute and the fiber product of $A$ and $X$ consists of a unique component.
Then for the integer
 \be \l{r} r=\vert  G_X\vert \vert Aut(G_X)\vert\ee the equality 
$B_1^{\circ r}=B_2^{\circ r}$
holds. 
\ec 
\pr 
Applying Theorem \ref{ite}, we can find $F_1,$ $F_2$ such that the diagrams  
\be 
\begin{CD}
S_X @> F_1 >> S_X \\
@VV H V @VV H V\\ 
C @>B_1 >> C \\
@VV X V @VV X V\\ 
\C\P^1 @>A  >> \C\P^1\,, 
\end{CD}\ \ \ \ \ \ \ \ \  \  
\begin{CD}
S_X @> F_2 >> S_X \\
@VV H V @VV H V\\ 
C @>B_2 >> C \\
@VV X V @VV X V\\ 
\C\P^1 @>A  >> \C\P^1\  
\end{CD}
\ee
commute. Furthermore, 
considering instead of the functions $A$, $B_1,$ $B_2,$ $F_1,$ $F_2$ their   $\vert Aut(G_X)\vert$-th iterates, we may assume that the corresponding automorphisms 
$\phi_i,$ $i=1,2$, of $S_X$ are the identical automorphisms, that is, that $F_i$, $i=1,2$, commute with $G_X$. 
Under this assumption, we must show that 
\be \l{eg2}  B_1^{\circ \vert G_X\vert}= B_2^{\circ \vert G_X\vert}.\ee
Since 
$$ F_2= g\circ  F_1$$
	for some  $g \in  G_X$, we have:
$$ F_2^{\circ \vert G_X\vert}=( g\circ F_1)^{\circ \vert G_X\vert}= g^{\circ \vert G_X\vert}\circ    F_1^{\circ \vert G_X\vert}= F_1^{\circ \vert G_X\vert},$$
implying \eqref{eg2}.

\bc \l{overq}
Let $E$ be an algebraic curve over $\C$, and $X:E\rightarrow \C\P^1$, $B:E\rightarrow E$, $A:\C\P^1\rightarrow \C\P^1$ 
dominant morphisms such that the diagram 
\be \l{xoma}
\begin{CD}
E @>B>> E \\
@VV X V @VV  X V\\ 
\C\P^1 @>A>> \C\P^1\
\end{CD}
\ee
commutes and the fiber product of $A$ and $X$ consists of a unique component.
 Assume that the curve $E$ and the morphisms $X,A$ are defined over some number field  $K$. 
 Then for the integer
 \be  r=\vert  G_X\vert \vert Aut(G_X)\vert\ee
the iterate $B^{\circ r}$ is  defined over $K.$
\ec

\pr 
It is clear that $B$ is defined over $\bar{\Q}$ and that for  any  $\gamma\in \Gal(\bar\Q/K)$ the function $^\gamma B$ satisfies \eqref{xoma} along with $B$. Thus, 
$$B^{\circ r}= (^\gamma B)^{\circ r},$$ by Corollary \ref{ite}.
Since  $$^\gamma \left(B^{\circ r}\right)= (^\gamma B)^{\circ r},$$ this implies that 
\be \l{esz} ^\gamma \left(B^{\circ r }\right)= B^{\circ r}\ee
  for any  $\gamma\in \Gal(\bar\Q/K)$, and hence,  $B^{\circ r}$ is  defined over $K.$ \qed

\end{subsection} 

\begin{subsection}{Proof of Theorem \ref{t2}} Without loss of generality, we may assume that  
$x_0$ is not $A$-preperiodic, since otherwise the theorem is obviously true. 
As in Section 3, we will consider the fiber products of $A^{\circ d},$ $d\geq 1,$ and $U$, but now regarding them as algebraic curves 
\be \f E_d:\, A^{\circ d}(x)-U(y)=0\ee in $(\C\P^1)^2.$
We observe that $\f E_d$, $d\geq 2,$ is the preimage of $\f E_{d-1}$ under the map  
$f:(\C\P^1)^2\rightarrow (\C\P^1)^2$ defined by 
$$f:\,(x,y)\rightarrow (A(x),y).$$ Let us denote by $\pi_x$  and $\pi_y$ the 
projection maps to $x$ and $y$ in $(\C\P^1)^2$.

Let $E_d,$ $d\geq 1,$ be a sequence of  
 irreducible components of $\f E_d$, $d\geq 1,$ such that   the diagram 
\be \l{moca}
\dots \begin{CD}
E_3@> f_3 >> E_2  @> f_2 >> E_1 @> \pi_{y,1} >> \C\P^1\\
@VV \pi_{x,3} V @VV \pi_{x,2} V @VV \pi_{x,1} V @VV U V\\
 \C\P^1 @>A >> \C\P^1  @> A  >> \C\P^1 @> A   >> \C\P^1\,,
\end{CD}
\ee where   $f_d$, $\pi_{x,d}$, $\pi_{y,d}$, $d\geq 1$, denotes the restriction of $f$, $\pi_x$, $\pi_y$ on the curve $E_d$, 
commutes. 
It is clear that  to prove Theorem \ref{t2}, it is enough to show that if the set $I$ consisting of $i\in \N$ such that the curve $E_i$ contains a $K$-point of the form $(x_0,y)$ is infinite,  then $I$ is a finite union of arithmetic progressions. 

For any $K$-point  $(x_0,y)$ of $ E_i$, the point $(A^{s}(x_0),y)$, where $s\leq i-1$, is the $K$-point  of $ E_{i-s}.$
Since by assumption  $x_0$ is not $A$-preperiodic, this implies that if $I$ is infinite, then every 
curve  $E_d,$ $d\geq 1,$ has infinitely many $K$-points, implying by the Faltings theorem that $g(E_d)\leq 1,$ $d\geq 1.$
Let us consider, along with diagram \eqref{moca}, the diagram 
\be \l{coma4}
\dots \begin{CD}
\t E_3@> F_3 >> \t E_2  @> F_2 >> \t E_1 \\
@VV \eta_3 V @VV \eta_2 V @VV \eta_1 V\\
E_3@> f_3 >> E_2  @> f_2 >> E_1 @>  \pi_{y,1} >> \C\P^1\\
@VV \pi_{x,3} V @VV \pi_{x,2} V @VV \pi_{x,1} V @VV U V\\
 \C\P^1 @>A >> \C\P^1  @> A  >> \C\P^1 @> A   >> \C\P^1,
\end{CD}
\ee
where $\eta_d: \t E_d\rightarrow E_d$, $d\geq 1,$ are desingularization maps, and $F_d:\t E_d \rightarrow \t E_{d-1}$, $d\geq 2,$ are 
holomorphic maps between compact Riemann surfaces. Since Corollary \ref{volk} and Corollary \ref{sum+} imply that there exists $d_0$ such that the diagram 
$$D[\t E_d,A,\pi_{x,d}\circ \eta_d,F_d], \ \ \ \ d\geq d_0,$$ is good, applying Theorem \ref{xyi}, we conclude that 
 there exist  $s_0\geq d_0$ and $l\geq 1$ such 
that for any $d\geq s_0$, the Riemann surfaces $\t E_d$ and $\t E_{d+l}$ are isomorphic and 
\be \l{lp} (\pi_{x,d}\circ \eta_d)=\pi_{x,d+l}\circ \eta_{d+l}\circ \t\alpha_{d},\ee for some isomorphism 
$\t\alpha_{d}:\t E_d\rightarrow \t E_{d+l}.$ 

Since $\t\alpha_{d}$
descends to an automorphism  $\alpha_{d}:E_d\rightarrow E_{d+l}$ that makes the diagram 
\be 
\begin{CD}
\t E_d @>\t \alpha_d>> \t E_{d+l} \\
@VV \eta_d V @VV  \eta_{d+l} V\\ 
E_d @>\alpha_d>> E_{d+l}\
\end{CD}
\ee
commutative,  it follows from \eqref{lp} that  for every $d\geq s_0$ the equality 
\be \l{nate} \pi_{x,d}=\pi_{x,d+l}\circ \alpha_d\ee holds and the diagram 
\be \l{avot}
\begin{CD}
E_d @>R_d>> E_{d} \\
@VV \pi_{x,d} V @VV \pi_{x,d}  V\\ 
\C\P^1   @>A^{\circ l} >> \C\P^1  \,,
\end{CD}
\ee
where $$R_d=f_{d+1} \circ \dots \circ f_{d+l-1}\circ f_{d+l}\circ \alpha_d,$$ 
commutes. 
By Corollary \ref{overq}, for every $d\geq s_0$  there exists $r$ such that $R_d^{\circ r}$  
is defined over $K$. Moreover, since $r$ is defined in terms of the monodromy group of $\pi_{x,d}$,  it follows from equality \eqref{nate} that 
for all $d\geq s_0$ from  the same class  by modulo $l$ we can take the same $r.$ 
Therefore, considering the least common multiple of the corresponding $r$ for all such classes, 
without loss of generality we may assume that  $R_d^{\circ r}$  
is defined over $K$ for all $d\geq s_0.$ 

Let us assume now that  $E_{i_0}$ contains a $K$-point of the form $(x_0,y)$ for $i_0\geq s_0.$ 
Since \eqref{avot} implies that for every $k\geq 1$  the equality 
$$A^{\circ{lk}}(x_0)=\pi_{x,i_0}\circ R_{i_0}^{\circ k}(x_0,y)$$ 
holds, setting $R=lr$, we have:
\begin{multline} 
A^{\circ (i_0+Rk)}(x_0)=A^{i_0}\circ A^{\circ Rk}(x_0)=A^{i_0}\circ\pi_{x,i_0}\circ R_{i_0}^{\circ rk}(x_0,y)=\\ =U\circ \pi_{y,1}\circ f_{2}\circ  f_{3}\circ... \circ f_{i_0}\circ  R_{i_0}^{\circ rk}(x_0,y).
\end{multline} 
Therefore, since $R_d^{\circ r}$ is defined over $K$, all the numbers $A^{\circ (i_0+Rk)}(x_0)$, $k\geq 1,$ belong to $U(K).$ 
This shows that the set of $i\in \N$ such that the curve $E_i$ contains a $K$-point of the form $(x_0,y)$ is a union of a finite set and a finite number of arithmetic progressions with denominator $R$. 

Finally, if $A$ is not a generalized Latt\`es map, then arguing as in  Theorem \ref{chai} we conclude that 
there exists $s_0\geq 0$ such that $\deg \pi_{x,d}=1$ for all $d\geq s_0$. Therefore, \eqref{nate} 
and \eqref{avot} hold for $l=1.$ Moreover, since $\deg \pi_{x,d}=1$, the map $R_d$ is defined over $K,$ 
 so that $A^{\circ (i_0+k)}(x_0)\in U(K)$ for all $k\geq 0.$ \qed

\end{subsection} 

\begin{subsection}{Example} 
In conclusion, we illustrate some of the constructions and results of this paper with the following example: 
$$A=144\,{\frac {z \left( z+3 \right) }{ \left( z-9 \right) ^{2}}},\ \ \ \ \ \ U=z^2,\ \ \ \ \ \ z_0=1, \ \ \ \ \ \ k=\Q.$$ 
The function $A$ is obtained from a one-parameter series 
 introduced in the paper \cite{jo} where the value of the parameter is equal to one. 
 It is shown in \cite{jo} that $$I=\{0,2\}\cup \{1+2m: m\geq 0\},$$
  so, by Theorem \ref{t2}, the function $A$ should be 
a generalized Latt\`es map. 
Specifically,  $A:\f O\rightarrow \f O$ is a minimal holomorphic map for the orbifold $\f O$ defined by 
the equalities 
$$\nu(0)=2, \ \ \ \ \nu(-3)=2.$$ Indeed, $$A^{-1}(0)=\{0,-3\},$$ and 
the multiplicity of $A$ at the points $0$ and $-3$ equals one, and hence, \eqref{uu} holds at $z=0$ and $z=-3.$  Moreover, 
$$A^{-1}(-3)=-9/7,$$ and 
 the multiplicity of $A$ at  $-9/7$ equals two, and hence, \eqref{uu} holds at $z=-9/7.$ 
On the other hand, for any point $z$ distinct from $0,-3,$ and $-9/7,$ equality  \eqref{uu} also holds, since 
for such a point,  $\nu(z)=1$ and $\nu(A(z))=1$.

To simplify formulas, let us consider instead of the functions $A$ and $U$ 
the functions 
$$\t A=\mu \circ A\circ \mu^{-1}, \ \ \ \  \t U=\mu \circ U,$$ where $$\mu=\frac{z}{z+3},$$ so that  
$$\t A=48\,{\frac {z}{ \left( 4\,z+3 \right) ^{2}}}, \ \ \ \ \t U={\frac {{z}^{2}}{{z}^{2}+3}}.$$ Then 
$\t A:\t{\f O}\rightarrow \t{\f O}$ is a minimal holomorphic map for the orbifold $\t{\f O}$ defined by the equalities 
$$\t\nu(0)=2, \ \ \ \ \t\nu(\infty)=2,$$ and the functions $F$ and $\theta_{\t{\f O}}$ from Proposition \ref{poiu}, which make the diagram 
\be 
\begin{CD}
 \C\P^1  @>F >>  \C\P^1\\
  @VV \theta_{\t{\f O}} V @VV \theta_{\t{\f O}} V\\
\C\P^1 @> \t A   >> \C\P^1
\end{CD}
\ee
commutative, have the form
$$F=4\,{\frac {\sqrt {3}z}{4\,{z}^{2}+3}}, \ \ \ \ \theta_{\t{\f O}}=z^2.$$

Further, the diagram
\be \l{xara0} 
\begin{CD}
 \C\P^1  @>F >>  \C\P^1\\
  @VV {(\theta_{\t{\f O}},V)} V @VV {(\theta_{\t{\f O}},V)} V\\
E  @> R >>  E @> \pi_y >>  \C\P^1\\
  @VV \pi_x V @VV \pi_x V @VV \t U V\\
\C\P^1 @> \t A   >> \C\P^1 @> \t A >> \C\P^1\,,
\end{CD}
\ee
where 
$$V=12\,{\frac {z}{4\,{z}^{2}-3}}$$ 
and the morphism $R:E\rightarrow E$ is defined by 
$$x\rightarrow \t A(x), \ \ \ \ y\rightarrow -48\,{\frac {\sqrt {3} \left(  \left( 4\,xy-3\,y \right) ^{2}+108
 \right)  \left( 4\,xy-3\,y \right) }{ \left(  \left( 4\,xy-3\,y
 \right) ^{2}-36 \right)  \left(  \left( 4\,xy-3\,y \right) ^{2}-324
 \right) }}\,,
$$ 
commutes. 
Finally, the function $\bar F$ Galois conjugated to $F$  satisfies $\bar F=- F,$ and 
the function  
$$ F^{\circ 2}= \bar F^{\circ 2}=16\,{\frac { \left( 4\,{z}^{2}+3 \right) z}{16\,{z}^{4}+88\,{z}^{2}+9}
}
$$ as well as the morphism $R^{\circ 2}$ have rational coefficients.

\vskip 0.2cm

\noindent {\bf Acknowledgments}. The  author is grateful for the hospitality and support he received at the Institute des Hautes \'Etudes Scientifiques,
where this work was begun, and at the Max-Planck-Institut
f\"ur Mathematik, where this
work was completed.

\end{subsection} 

\end{section}

\end{document}